\documentclass[11pt]{amsart}
\usepackage{amsmath}
\usepackage{amscd}
\usepackage{graphics}
\usepackage{latexsym}

\oddsidemargin .in
\theoremstyle{plain}
\newtheorem{Thm}{Theorem}
\newtheorem{Lem}[Thm]{Lemma}

\newtheorem{Prop}[Thm]{Proposition}
\newtheorem{Rm}{Remark}
\newtheorem{Def}{Definition}

\theoremstyle{remark}
\def\bL{{\bf L}}
\def\bE{{\bf E}}
\def\bV{{\bf V}}
\def\bW{{\bf W}}
\def\bK{{\bf K}}

\def\C{{\mathbb C}}
\def\A{{\mathbb A}}

\def\E{{\mathbb E}}
\def\R{{\mathbb R}}
\def\L{{\mathbb L}}
\def\P{{\mathbb P}}

\def\Z{{\mathbb Z}}
\def\W{{\mathbb W}}

\def\O{{\mathbb O}}
\def\H{{\mathbb H}}
\def\V{{\mathbb V}}

\def\CM{\mathcal M}
\def\CG{\mathcal G}

\def\CA{\mathcal A}

\def\CP{\mathcal P}
\def\CS{\mathcal S}
\def\CZ{\mathcal Z}

\def\CL{\mathcal L}

\def\s2x{\hbox{$S^2 \times S^2$}}

\def \Di {D\!\!\!\!/\,}

    \def\sqr#1#2{{\vcenter{\hrule height.#2pt
            \hbox{\vrule width.#2pt height#1pt \kern#1pt
            \vrule width.#2pt}\hrule height.#2pt}}}
    \def\square{\mathchoice\sqr67\sqr67\sqr{2.1}6\sqr{1.5}6}

\def\qed{~\hfill$\square$}

\begin{document}

\title[]{Calibrated Manifolds and Gauge theory  }
\author{Selman Akbulut and Sema Salur}
\thanks{First named author is partially supported by NSF grant DMS 9971440}
\address{Department  of Mathematics, Michigan State University,  MI, 48824}
\email{akbulut@math.msu.edu }
\address {Department. of Mathematics, Northwestern University, IL, 60208 }
\email{salur@math.northwestern.edu }
\date{\today}

\begin{abstract}
By a theorem of Mclean, the deformation space of an associative
submanifold $Y$ of an integrable $G_2$-manifold $(M,\varphi)$ can
be identified with the kernel of a Dirac operator
$\Di:\Omega^{0}(\nu)\to \Omega^{0}(\nu)$ on the normal bundle
$\nu$ of $Y$. Here, we generalize this to the non-integrable case,
and also show that the deformation space becomes smooth after
perturbing it by natural parameters, which corresponds to moving
$Y$ through `pseudo-associative' submanifolds. Infinitesimally,
this corresponds to twisting the Dirac operator $\Di\mapsto
\Di_{A}$  with connections $A$ of $\nu$. Furthermore, the normal
bundles of the associative submanifolds with $Spin^c$ structure
have natural complex structures, which helps us to relate their
deformations to Seiberg-Witten type equations.

If we consider $G_2$ manifolds with $2$-plane fields $(M,\varphi,
\Lambda)$ (they always exist) we can split the tangent space $TM$
as a direct sum of an associative 3-plane bundle and a complex
$4$-plane bundle. This allows us to define (almost) $\Lambda$-associative
submanifolds of $M$, whose deformation equations, when perturbed,
reduce  to Seiberg-Witten equations, hence we can assign local
invariants to these submanifolds. Using this we can assign
an invariant to $(M,\varphi, \Lambda)$.  These Seiberg-Witten
equations on the submanifolds are  restrictions of global
equations on $M$.
 We also discuss similar results for the Cayley
submanifolds of a $Spin(7)$ manifold.
\end{abstract}

\maketitle

\setcounter{section}{-1}

\section{ Introduction}

\vspace{.1in}

We first study deformations of associative submanifolds $Y^3$ of a
$G_2$ manifold $(M^7,\varphi )$, where $\varphi \in \Omega^{3}(M)$
is  the $G_2$ structure.  We prove a generalized version of the
McLean's theorem where integrability condition of the underlying
$G_2$ structure is not necessary. This deformation space might be
singular, but  by perturbing it with some natural parameters it
can be made smooth. This amounts to deforming $Y$ through the
associatives in $(M,\varphi)$ with varying $\varphi $,  or
alternatively deforming $Y$ through the pseudo-associative
submanifolds  ($Y$'s whose tangent planes become associative after
rotating by a generic element of the gauge group of $TM$).
Infinitesimally, these perturbed deformations correspond to the
kernel of the twisted Dirac operator  $\Di_{A}:\Omega^{0}(\nu)\to
\Omega^{0}(\nu)$, twisted by some connection $A$ in $\nu(Y)$.

\vspace{.1in}

The associative submanifolds with $Spin^c$ structures in
$(M,\varphi)$ are useful objects to study, because their normal
bundles have natural complex structures. Also we can view $ (M,
\varphi)$  as an analog  of a symplectic manifold,  and view a
non-vanishing $2$-plane field $\Lambda$ on $M$ as an analog of a
complex structure taming $\varphi$. Note that $2$-plane fields are
stronger versions of $Spin^{c}$ structures on $M^7$, and they
always exist by \cite{t}. The  data $(M^7,\varphi, \Lambda)$
determines an interesting  splitting of the tangent bundle $TM=
\bE\oplus \bV$, where $\bE $ is the bundle of associative
$3$-planes, and $\bV$ is the complementary $4$-plane bundle with a
complex structure, which is a spinor bundle of $\bE$. Then the
integral submanifolds $Y^3$ of $\bE$, which we call {\it
$\Lambda$-associative submanifolds},  can be viewed as analogues of
J-holomorphic curves; because their normal bundles come with an almost complex structure. Even if they may not always exist, their  perturbed
versions, i.e. {\it almost $\Lambda$-associative submanifolds}, always do.
Almost $\Lambda$-associative submanifolds are the transverse sections of
the bundle ${\bf V}\to M$. We can deform such $Y$ by using the
connections in the determinant line bundle of  $\nu(Y)$ and get a
smooth deformation space, which is described by the  twisted Dirac
equation. Then by constraining this new variable with another natural equation
we arrive to  Seiberg-Witten type equations for $Y$. So we can
assign an integer to $Y$, which is invariant under small isotopies through almost $\Lambda$-assocative submanifolds.

\vspace{.1in}

In fact it turns out that $(M^7,\varphi, \Lambda)$ gives a finer
splitting $TM=\bar{\bE}\oplus \xi$, where $\bar{\bE} $ is a
$6$-plane bundle with a complex structure, and $\xi$ is a real
line bundle. In a way this structure of $(M, \varphi)$ mimics the
structure of (Calabi-Yau)$ \times S^1$ manifolds, and by
`rotating' $\xi$ inside of $TM$ we get a new insight for so-called
``Mirror manifolds" which is investigated  in \cite{as1}.

\vspace{.1in}

There is a similar process for the deformations of Cayley
submanifolds $X^4\subset N^8$  of a Spin(7) manifold $(N^8,\Psi)$,
which we discuss at the end. So in a way $\Lambda$-associative  (or
Cayley) manifolds in a $G_2$ (or $Spin(7)$) manifold, behave much
like higher dimensional analogue of holomorphic curves in a
Calabi-Yau manifold.

\vspace{.05in}

We would like to thank MSRI, IAS, Princeton and Harvard
Universities for providing a stimulating  environment where this
paper is written, and we thank R. Kirby and G. Tian  for
continuous encouragement. The first named author thanks to R.
Bryant and C. Taubes for stimulating discussions and useful
suggestions.

\newpage

\section{ Preliminaries}

\vspace{.1in}

Here we first review basic properties of the manifolds with
special holonomy (most material can be found in \cite{b2},
\cite{b3}, \cite{h}, \cite{hl}), and then proceed to prove some
new results. Recall that the set of octonions $\O=\H\oplus l \H=
\R^8$  is an $8$-dimensional division algebra generated by $<1, i,
j,  k, l, li ,lj, lk> $. On  the set of the imaginary octonions
$im \O =\R^7$ we have the cross product operation $\times:
\R^7\times \R^7 \to \R^7$, defined  by $u\times v=im(\bar{v}.u)$.
The  exceptional Lie group $G_{2}$ can be defined as the linear
automorphisms of  $im \O$ preserving this cross product operation,
$G_2=Aut(\R^7, \times)$. There is also another  useful description
in terms of the orthogonal $3$-frames in $\R^7$:
\begin{equation}
G_{2} =\{ (u_{1},u_{2},u_{3})\in (im \O)^3\;|\: <u_{i},u_{j}>=\delta_{ij},  \;
<u_{1} \times u_{2},u_{3}>=0 \;  \}
\end{equation}

Alternatively,  $G_2$ can be defined as the  subgroup of the
linear group $GL(7,\R)$ which  fixes  a particular $3$-form
$\varphi_{0} \in \Omega^{3}(\R ^{7})$.  Denote
$e^{ijk}=dx^{i}\wedge dx^{j} \wedge dx^{k}\in \Omega^{3}(\R^7)$,
then

$$ G_{2}=\{ A \in GL(7,\R) \; | \; A^{*} \varphi_{0} =\varphi_{0}\; \} $$
\vspace{-0.15in}
\begin{equation}
\varphi _{0} =e^{123}+e^{145}+e^{167}+e^{246}-e^{257}-e^{347}-e^{356}
\end{equation}

{\Def A smooth $7$-manifold $M^7$ has a {\it $G_{2}$ structure} if its tangent frame bundle reduces to a $G_{2}$ bundle. Equivalently,
 $M^7$ has a {\it $G_{2}$ structure} if there is  a 3-form $\varphi \in \Omega^{3}(M)$  such that  at each $x\in  M$  the pair $ (T_{x}(M), \varphi (x) )$ is  isomorphic to
 $(T_{0}( \R^{7}), \varphi_{0})$}.

 \vspace{.1in}

Here are some useful properties, discussed more fully in \cite{b2}:  Any  $G_{2}$ structure $\varphi$ on $M^7$ gives an orientation $\mu \in \Omega^{7}(M)$ on $M$,   and this $\mu$  determines  a metric $g= \langle \;,\;\rangle$ on $M$, and a cross product structure $\times$  on its tangent bundle of $M$ as follows: Let  $i_{v}$ denote the interior product with a vector $v$ then
\begin{equation}
\langle u,v \rangle=[ i_{u}(\varphi ) \wedge i_{v}(\varphi )\wedge \varphi  ]/6\mu
\end{equation}
\begin{equation}
\varphi (u,v,w)=\langle u\times v,w \rangle
\end{equation}

\noindent To emphasize the dependency on $\varphi $ sometimes $g$
is denoted by $g_{\varphi}$. In particular, the $14$-dimensional
Lie group $G_{2}$  imbeds into $SO(7)$  subgroup of $GL(7, \R)$.
Note that because of the way we defined $G_2=G_2^{\varphi_{0}}$,
this imbedding is determined  by $\varphi_{0}$.

 \vspace{.05in}

Since $GL(7,\R)$ acts on $\Lambda^{3}(\R^7)$ with  stabilizer
$G_{2}$, its orbit $\Lambda^{3}_{+}(\R^7) $ is open for dimension
reasons, so the choice of $\varphi_{0}$ in  the above definition
is generic (in fact it has two orbits  containing $\pm
\varphi_{0}$). $G_2$  has many copies $G_{2}^{\varphi}$  inside
$GL(7, \R)$,  which are all conjugate to each other, since $G_2$
has only one $7$ dimensional representation.  Hence the space of
$G_2$ structures on $M^7$ are identified with the sections of the
bundle:
\begin{equation}
 \R\P^7 \simeq  GL(7,\R) /G_{2}\to \Lambda ^{3}_{+}(M) \longrightarrow M
 \end{equation}
which are called {\it the positive $3$-forms}, these are the set of $3$-forms $\Omega^{3}_{+}(M)$  that can  be identified  pointwise by $\varphi_{0}$.
Each $G_{2}^{\varphi }$ imbeds into  a conjugate  of one standard copy $SO(7)\subset GL(7, \R)$. The
space of $G_2$ structures $ \varphi $ on $M$, which induce the same metric on $M$, that is all $\varphi $'s for which the corresponding $G_{2}^{\varphi}$ lies in the standard $SO(7)$,  are  the sections of the bundle (whose fiber is the orbit of $\varphi_{0}$ under $SO(7)$):
\begin{equation}
 \R\P^7 = SO(7)/G_{2}\to  \tilde{\Lambda}^{3}_{+}(M) \longrightarrow M
 \end{equation}
 which we will denote by $\tilde{\Omega}^{3}_{+}(M)$.
 The set of smooth $7$-manifolds with $G_{2}$-structures coincides with  the set of $7$-manifolds with spin structure, though this correspondence is not $1-1$. This is because $Spin(7)$ acts on $S^{7}$ with stabilizer $G_{2}$ inducing  the fibrations
 \begin{equation*}
G_{2}\to Spin(7)\to S^{7}\to BG_{2}\to BSpin(7)
\end{equation*}
and so there is no obstruction to lifting maps $M^7\to BSpin(7)$ to $BG_{2}$,  and there are many liftings.
Cotangent frame bundle $\CP^{*}(M)\to M$ of a manifold with $G_{2}$ structure $(M,\varphi)$  can be expressed as $\CP^{*}(M)=\cup_{x\in M} \;\CP^{*}_{x}(M)$, where
each fiber is:
\begin{equation*}
\CP^{*}_{x}(M)=\{ u\in Hom (T_{x}(M), {\R}^{7})\;|\; u^{*}(\varphi_{0})=\varphi (x)\;\}
\end{equation*}

 \vspace{.1in}

Throughout this paper we will denote the cotangent frame bundle by
$\CP^{*}(M)\to M$ and its adapted frame bundle by $\CP(M)$. They
can be $G_{2}$ or $SO(7)$ frame bundles; to emphasize it sometimes
we will specify them by the notations  $\CP_{SO(7)}(M)$ or
$\CP_{G_{2}}(M)$. Also  we will denote the sections of a bundle
$\xi \to Y$ by $\Omega^{0}(Y,\xi)$ or simply by $\Omega^{0}(\xi)$,
and the bundle valued $p$-forms by
$\Omega^{p}(\xi)=\Omega^{0}(\Lambda^{p}T^{*}Y\otimes \xi)$, and
the {\it sphere bundle} of $\xi$ by $S(\xi)$. There is a notion of
a $G_2$ structure $\varphi $  on $M^7$ being {\it integrable},
which corresponds to $\varphi$ being an harmonic form:

%\vspace{.05in}

{\Def A manifold with $G_{2}$ structure $(M,\varphi)$  is called a {\it $G_{2}$ manifold} if   the holonomy group of the Levi-Civita connection (of the metric $g_{\varphi }$) lies inside of  $G_2 $. Equivalently  $(M,\varphi)$ is a $G_{2}$ manifold if $\varphi $ is parallel with respect to the metric $g_{\varphi }$ i.e.
$\nabla_{g_{\varphi }}(\varphi)=0$; this condition is equivalent to $d\varphi=0 =\;d(*_{g_{\varphi}}\varphi) $.}

\vspace{.1in}

In short one can define a $G_{2}$ manifold to be any Riemannian manifold $(M^{7},g)$ whose holonomy group is contained in $G_{2}$,  then $\varphi $ and the cross product $\times$ come  as a consequence. It turns out that  the condition $\varphi$ being harmonic is equivalent to the condition that at each point $x_{0}\in M$ there is a chart  $(U,x_{0}) \to (\R^{7},0)$ on which $\varphi $ equals to $\varphi_{0}$ up to second order term, i.e.  on the image of $U$
\begin{equation}
\varphi (x)=\varphi_{0} + O(|x|^2)
\end{equation}

\vspace{.1in}

{\Rm
For example if  $(X^6,\omega, \Omega)$ is a complex 3-dimensional Calabi-Yau manifold with
K\"{a}hler form $\omega$, and a nowhere vanishing holomorphic
3-form $\Omega$, then $X\times S^1$ has holonomy group
$SU(3)\subset G_2$, hence  is a $G_2$ manifold. In this case
\begin{equation}
\varphi= \mbox{Re} \; \Omega + \omega \wedge dt.
\end{equation}}

%\vspace{.05in}

{\Def  Let $(M, \varphi )$ be a manifold with a $G_2$ structure. A
4-dimensional submanifold $X\subset M$ is called  an {\em
co-associative } if $\varphi|_X=0$. A 3-dimensional submanifold
$Y\subset M$ is called  an {\em associative} if $\varphi|_Y\equiv
vol(Y)$; this condition is equivalent to $\chi|_Y\equiv 0$,  where
$\chi \in \Omega^{3}(M, TM)$ is the  tangent bundle valued 3-form
defined   by  the identity: }
\begin{equation}
\langle \chi (u,v,w) , z \rangle=*\varphi  (u,v,w,z)
\end{equation}
The equivalence of these  conditions follows from  the `associator equality' of  \cite{hl}
\begin{equation}
\varphi  (u,v,w)^2 + |\chi (u,v,w)|^2/4= |u\wedge v\wedge w|^2
\end{equation}

 In general, if $\{e^1,e^2,..,e^7\}$ is any orthonormal coframe on $(M, \varphi)$, then the expression  (2) for $\varphi$ hold on a chart.
  By calculation $*\varphi $, and using (9) we can  calculate the expression of $\chi $ (note the error in the the second term of $6$th line of the corresponding  formula (5.4) of \cite{m}):
\begin{equation} *\varphi  =e^{4567}+e^{2367}+e^{2345}+e^{1357}-e^{1346}-e^{1256}-e^{1247}
\end{equation}
\begin{eqnarray*}
\chi  \hspace{-.1in}  &=&\;(e^{256}+ e^{247}+e^{346} -e^{357} ) \;e_{1}\\
            &&  + \; ( -e^{156}- e^{147}-e^{345} -e^{367} ) \;e_{2} \\
             &&+ \; (e^{245}+ e^{267}-e^{146} +e^{157} ) \;e_{3} \\
             &&+ \;(-e^{567}+ e^{127}+e^{136} -e^{235} ) \;e_{4} \\
            && + \;(e^{126}+ e^{467}-e^{137} +e^{234} ) \;e_{5} \\
            &&+ \;(-e^{457}- e^{125}-e^{134} -e^{237} ) \;e_{6}  \\
             &&+ \;(e^{135} -e^{124} +e^{456}+e^{236} ) \;e_{7}
\end{eqnarray*}

\noindent  Also $\chi$ can be expressed in terms of cross product operation (c.f. \cite{h}, \cite{hl}, \cite{k}):
\begin{equation}
\chi(u,v,w)= -u\times (v\times w)-\langle u,v\rangle w +\langle u,w\rangle v
\end{equation}

 \vspace{.05in}

 When $d\varphi=0$, the associative submanifolds are volume minimizing submanifolds of $M$ (calibrated by $\varphi $).  Even in  the general case of  a manifold with a $G_2$ structure  $(M,\varphi )$, the form $\chi$ imposes an interesting structure near associative submanifolds:

  \vspace{.1in}

Notice  (9) implies that,  $\chi $ maps every oriented $3$-plane in $T_{x}(M)$ to the orthogonal subspace $T_{x}(M)^{\perp}$, so  if we choose local coordinates  $(x_1,..., x_7)$ for $M^7$ we get
\begin{equation}
\chi =\sum  a^{\alpha }_{J } \; dx^{J}
\otimes \frac{\partial}{\partial x_{\alpha}}
\end{equation}
where $dx^{J}= dx^{i}\wedge dx^{j}\wedge dx^{k}$,  and  the summation is taken over the multi-index  $J=\{i,j,k\}$  and $\alpha$ such that $\alpha\notin J$.  So if $ Y\subset M $ is given by $ (x_1,x_2,x_3 )$ coordinates, then locally the condition $Y$ to be associative is given by the equations:
\begin{equation}
a^{\alpha }_{123 }=0
\end{equation}
From (9) it is easy to calculate $a^{\alpha }_{ijk}= *\varphi_{ijks} g^{s \alpha }$, where
$g^{-1}=(g^{ij })$ is the inverse of the metric  $g=(g_{ij})$, and of course the metric $g $ can be expressed in terms of $\varphi $. By evaluating $\chi$ on the orientation form of $Y$ we get a normal vector field so:

{\Lem  To any $3$-dimensional submanifold $Y^3\subset (M,\varphi)$,  $\chi$ associates a normal vector field, which vanishes when $Y$ is associative.}

\vspace{.1in}

Hence $\chi $ defines an interesting flow on $3$ dimensional submanifolds of $(M,\varphi )$, fixing associative submanifolds. On the associative submanifolds with a $Spin^c$ structure, $\chi $ rotates their normal bundles and imposes a complex structure on them:

%\vspace{.05in}

{\Lem  To any associative manifold $Y^3\subset (M,\varphi)$ with a non-vanishing oriented $2$-plane field,  $\chi$ defines an almost complex structure on  its normal bundle $\nu(Y)$ (notice that in particular any coassociative submanifold $ X\subset M$ has an almost complex structure if its normal bundle has a non-vanishing section).}

\proof Let $L\subset \R^7$ be an associative $3$-plane, that is $\varphi |_{L}=vol(L)$. Then to every pair of orthonormal vectors $\{u,v\}\subset L$, the form $\chi$ defines a complex structure on the orthogonal $4$-plane $L^{\perp}$, as follows:   Define $j: L^{\perp} \to L^{\perp}$ by
\begin{equation}
j(X)=\chi(u,v,X)
\end{equation}
This is well defined i.e. $j(X)\in L^{\perp}$, because when $ w\in L$ we have:
$$ <\chi(u,v,X),w>=*\varphi (u,v,X,w)=-*\varphi (u,v,w,X)=<\chi(u,v,w),X>=0$$

\noindent Also $j^{2}(X)=j(\chi(u,v,X))=\chi(u,v,\chi(u,v,X))=-X$. We can check the last equality by taking an orthonormal basis $\{ X_{j}\}\subset L^{\perp}$ and calculating
\begin{eqnarray*}
<\chi(u,v,\chi(u,v,X_{i})),X_{j}>&=&*\varphi (u,v,\chi(u,v,X_{i}),X_{j})=\\
-*\varphi (u,v,X_{j},\chi(u,v,X_{i})) &=&- <\chi(u,v,X_{j}),\chi(u,v,X_{i})>=-\delta_{ij}
\end{eqnarray*}

The last equality holds since the map $j$ is orthogonal, and the orthogonality  can be seen by polarizing the associator equality (10), and by noticing $\varphi (u,v,X_i)=0$. Observe that the map $j$ only depends on the oriented $2$-plane $l=<u,v>$  generated by $\{u,v\}$. So the result follows. \qed

\vspace{.1in}

In fact, for any unit vector field $\xi$ on an associative  $Y$ (i.e. a $Spin^c$ structure)
 defines a complex structure $J_{\xi}: \nu(Y) \to \nu(Y) $ by $J_{\xi}(z)=z\times \xi$, and the complex structure defined in Lemma 2 corresponds to $J_{u\times v}$, because from (12):
 $$ \chi(u,v,z)=\chi(z,u,v)=
-z\times (u\times v)-\langle z, u\rangle v + \langle z, v\rangle u=J_{v\times u}(z).$$

\vspace{.05in}

\noindent Also recall that the complex structures on any $SO(4)$ bundle such as $\nu \to Y$ are given by the unit sections of the associated $SO(3)$  bundle $\lambda_{+}(\nu)\to Y$, which is induced by the left reductions $SO(4)=(SU(2)\times SU(2))/\Z_2\to SU(2)/\Z_2= SO(3)$.

\vspace{.1in}

{\Def A Riemannian 8-manifold $(N^{8} , g)$ is
called a  {\em Spin(7) manifold} if the holonomy group of its
Levi-Civita connection lies in $Spin(7)\subset GL(8,\R)$.}

\vspace{.05in}

Equivalently  a $ Spin(7) $ manifold $(N,\Psi) $ is  a Riemannian $8$-manifold with a triple  cross product $\times $ on its tangent bundle, and a harmonic $4$-form $\Psi \in \Omega^{4}(N)$  with
$$ \Psi (u,v,w,z)=g(u \times v \times w,z)$$

 It is easily checked that  if  $(M,\varphi )$ is a $G_2$ manifold, then $(M\times S^1, \Psi )$ is a $Spin(7)$ manifold where $\Psi= \varphi\wedge dt - * \varphi$.

{\Def  A 4-dimensional submanifold $X$ of a $Spin(7)$ manifold $(N,\Psi)$
is called {\em Cayley} if $\Psi|_X\equiv vol(X)$. This is equivalent to $\tau|_X\equiv 0$ where
$\tau\in \Omega^4(N,E)$ is a certain vector-bundle valued 4-form defined by the ``four-fold
cross product'' of the imaginary octonions $\tau (v_1,v_2,v_3,v_4)=v_1 \times v_2 \times v_3  \times v_4$ (see \cite{m}, \cite{hl}). }
 \vspace{.1in}

\vspace{.1in}

\section{ Grassmann Bundles }

\vspace{.1in}

Let $G(3,7)$ be the Grassmann manifold  of oriented $3$-planes in
$\R^7$. Let $M^7$ be an oriented smooth $7$-manifold, and let
$\tilde{M}\to M$ be the bundle oriented $3$-planes in $TM$, which
is defined by the identification  $[p,L]=[pg,g^{-1}L] \in
\tilde{M}$:

\begin{equation}
 \tilde{M}= \CP_{SO(7)}(M)\times_{SO(7)} G(3,7) \to M.
\end{equation}

\noindent This is  just the bundle $\tilde{M}=\CP_{SO(7)}(M)/
SO(3)\times SO(4) \to \CP_{SO(7)}(M)/SO(7) = M$.  Let $\xi \to
G(3,7)$ be the universal $\R^3$ bundle, and $\nu =\xi^{\perp}\to
G(3,7)$ be the dual $\R^{4}$ bundle. Therefore, $Hom(\xi,
\nu)=\xi^{*}\otimes \nu \longrightarrow G(3,7)$ is the tangent
bundle  $TG(3,7)$. $\xi$, $\nu$ extend fiberwise to give bundles
$\Xi \to  \tilde{M}$, $\V \to \tilde{M}$ respectively, and let
$\Xi^*$ be the dual of $\Xi$. Notice   that  $ Hom (\Xi, \V)=
\Xi^{*}\otimes \V \to \tilde{M}\; $ is the bundle  of vertical
vectors  $T^{v}(\tilde{M}) $ of $T(\tilde{M}) \to M$, i.e. the
tangents to the fibers of  $\pi: \tilde{M}\to M$, hence

\begin{equation}
T\tilde{M}\cong  T^{v}(\tilde{M}) \oplus  \pi^{*} TM   = (
\Xi^{*}\otimes \V ) \oplus  \Xi \oplus \V.
\end{equation}

\noindent That is, $T\tilde{M} $ is the  vector bundle associated
to principal $SO(3)\times SO(4)$ bundle $\CP_{SO(7)}\to \tilde{M}$
by the obvious representation of $SO(3)\times SO(4) $ to
$(\R^3)^{*}\otimes \R^{4} + \R^{3} + \R^{4}$.  The identification
(17) is defined up to gauge automorphisms of  bundles $\Xi$ and
$\V$.

\vspace{.1in}

Note that the bundle $\V=\Xi^{\perp}$ depends on the metric, and  hence it depends on $\varphi$ when metric is induced from a $G_2$ structure  $(M,\varphi )$. To emphasize this fact we can denote it by $\V_{\varphi} \to \tilde{M}$. But when we are considering $G_2$ structures coming from $G_2$ subgroups of a fixed copy of $SO(7)\subset GL(7, \R)$, they induce the same metric and so this distinction is not necessary.

\vspace{.05 in}

Let  $\CP(\V)\to \tilde{M}$ be the $SO(4)$ frame bundle of the vector bundle $\V$,  identify $\R^4$ with  the quaternions $ \H $, and identify $SU(2)$ with the unit quaternions $Sp(1)=S^3$. Recall  that $SO(4)$ is the equivalence classes of pairs $[\;q,\lambda \;]$ of unit quaternions
$$SO(4)=(SU(2)\times SU(2))/\Z_{2}$$ Hence   $\V\to \tilde{M} $  is the associated vector bundle to $\CP(\V)$ via  the $SO(4)$ representation
\begin{equation}
x \mapsto  qx\lambda^{-1}
\end{equation}
 There is  a pair of $\R^3=im(\H)$ bundles over $\tilde{M}$ corresponding to the left and right $SO(3)$ reductions of $SO(4)$, which are given by the $SO(3)$ representations
\begin{equation}
\begin{array}{lcc}
 \lambda_{+}(\V)\;:  &x  \mapsto qx\ q^{-1} &   \\
\lambda_{-}(\V)  \;: &y \mapsto  \lambda y  \lambda ^{-1} &
 \end{array}
 \end{equation}
The map $x\otimes y\mapsto xy$  gives  actions $\lambda_{+}(\V)\otimes \V \to \V  \; \hbox {and}\; \V\otimes \lambda_{-}(\V) \to \V$;   by combining we can think of them as one conjugation  action
\begin{equation}
(\lambda_{+}(\V)\otimes \lambda_{-}(\V))\otimes  \V \to \V
\end{equation}

\vspace{.1in}

If  the $SO(4)$ bundle $\CP(\V) \to \tilde{M}$ lifts to a $Spin(4)= SU(2)\times SU(2)$ bundle (locally it does),  we get two additional bundles  over $\tilde{M}$
 \begin{equation}
\begin{array}{lcc}
\;\CS\;: & y \mapsto  q y  & \\
\;\E\;: & \; \; y \mapsto  y  \lambda ^{-1} & \\
 \end{array}
 \end{equation}

\noindent They identify $\V$  as a  tensor product of two
quaternionic line bundles  $\V= \CS \otimes_{\H} \E$.  In
particular, $\lambda_{+}(\V)=ad(\CS)$ and $\lambda_{-}(\V)=
ad(\E)$, i.e. they are the $SO(3)$ reductions of the $SU(2)$
bundles $\CS$ and  $\E$.   Also there is  a multiplication map $
\CS \otimes \E \to \V$. Recall the  identifications:
$\Lambda^{2}(\V)= \Lambda^{2}_{+}(\V) \oplus \Lambda^{2}_{-}(\V)=
\lambda_{-} (\V)\oplus \lambda_{+}(\V) =\lambda(\V)=gl(\V)=ad
(\V)$.

\vspace{.1in}

\subsection{ Associative Grassmann Bundles}

$\;$

\vspace{.1in}

  Now consider  the {\it Grassmannian of associative $3$-planes} $G^{\varphi}(3,7) $  in $\R^7 $, consisting of elements $L\in  G(3,7)$ with the property $\varphi_{0}|_{L}=vol(L)$ (or equivalently $\chi_{0}|_{L}=0$). $G_2$ acts on $G^{\varphi}(3,7) $  transitively  with  the stabilizer $SO(4)$, so it gives the identification $G^{\varphi}(3,7)=G_{2}/SO(4) $.  If we identify  the imaginary octonions by $\R^7 =\mbox{Im} (\O )\cong im (\H)\oplus \H $,  then the action of the  subgroup $SO(4)\subset G_{2}$ on $\R^7$ is
\begin{equation}
\left(
\begin{array}{cc}
 \rho(A)  &  0 \\
   0  &   A  \\
\end{array}
\right)
\end{equation}
where $\rho :  SO(4) = (SU(2)\times SU(2))/ \Z_{2} \to SO(3)$ is the projection of the first factor (\cite{hl}), that is  for $[q,\lambda]\in SO(4)$ the action is given by  $(x,y)\mapsto (qxq^{-1}, qy\lambda^{-1})$.  So the action of $SO(4)$ on the $3$-plane $L=im (\H)$ is determined by its action on $L^{\perp}$. Now let  $M^7$  be a $G_2$ manifold.  Similar to the construction before,   we can construct  the bundle of associative Grassmannians over $M$ (which is a submanifold of $\tilde{M}$):
\begin{equation}
\tilde{M}_{\varphi} = \CP_{G_{2}}(M) \times _{G_{2}}G^{\varphi}(3,7) \to M
\end{equation}
which is just the quotient bundle $\tilde {M}_{\varphi}=\CP_{G_{2}}(M)/SO(4)\longrightarrow \CP_{G_{2}}(M)/G_{2} =M$. As in the previous section,  the restriction  of the universal bundles $\xi, \;\nu=\xi^{\perp} \to G^{\varphi}(3,7)$ induce $3$  and $4$ plane bundles $\Xi\to\tilde{ M}_{\varphi}$ and $\V\to \tilde{M}_{\varphi}$ (by restricting from  $ \tilde{M}$).   Also
\begin{equation}
T\tilde{M}_{\varphi }\cong T^{v } (\tilde{M}_{\varphi }) \oplus  \Xi \oplus \V
\end{equation}

From (22) we see that in  the  associative case,  we have an important identification: $\Xi= \lambda_{+}(\V)$  (as bundles over $\tilde{M}_{\varphi}$), and the dual of the  action $\lambda_{+}(\V)\otimes \V \to  \V$ gives a Clifford multiplication:
\begin{equation}
\Xi^{*} \otimes \V \to \V
\end{equation}

 In fact this is just the map induced from the cross product operation \cite{as2}.
     Recall that $T^{v}(\tilde{M})=  \Xi^{*}\otimes \V \to \tilde{M} $ is the  subbundle of  vertical  vectors of  $T (\tilde {M}) \to M$.   The total space $E(\nu_{\varphi} )$ of the normal bundle of the imbedding $ \tilde {M}_{\varphi } \subset \tilde{M}$  should be thought of an open tubular neighborhood of $ \tilde {M}_{\varphi }$ in $\tilde{M}$, and it has a nice description:

 \vspace{.05in}

{\Lem  (\cite{m})  Normal bundle $\nu_{\varphi}$ of $ \tilde {M}_{\varphi } \subset \tilde{M}$  is isomorphic to $\V$,  and the bundle of vertical vectors $T^{v} (\tilde{M}_{\varphi })$  is the kernel of the Clifford multiplication $ c: \Xi ^{*} \otimes \V \to \V$.
We have $ \; T^{v}(\tilde{M})| _{{\tilde{M}_{\varphi}}} =T^{v }(\tilde{M}_{\varphi })\oplus \nu _{\varphi }$,  and  the  following exact sequence over $\tilde{M}_{\varphi }$
  \begin{equation*}
T^{v} (\tilde{M}_{\varphi })\to \Xi^{*}\otimes \V|_{{\tilde{M}_{\varphi}}} \stackrel{c}{\longrightarrow} \V|_{\tilde{M}_{\varphi }} \to 0
\end{equation*}
\noindent Hence  the  quotient bundle, $ T ^{v} (\tilde{M}) / T^{v} (\tilde{M}_{\varphi })$  is isomorphic to $\V$.
\proof
This is because  the Lie algebra inclusion $ g_2\subset so(7) $ is given by
\[
\left(
\begin{array}{ccc}
a  &  \beta    \\
-\beta^t  &\rho ( a )    \\
 \end{array}
\right)
\]
\noindent  where $a\in so(4)$  is $y\mapsto qy-y\lambda$, and  $\rho(a)\in so(3)$ is
 $x\mapsto qx-xq$. So the tangent space inclusion of  $G_{2}/SO(4)\subset SO(7)/SO(4)\times SO(3)$ is given by the   matrix $\beta \in (im \H)^{*}\otimes \H $. Therefore,  if we write $\beta$ as column vectors  of three queternions $\beta=(\beta_1, \beta_2, \beta_3)=i^{*}\otimes \beta_1 +j^{*}\otimes \beta_2+ k^{*}\otimes \beta_{3}$, then $\beta_{1}i+\beta_{2}j+\beta_{3}k=0$ (\cite{m},  \cite{mc}). \qed

\vspace{.1in}

The reader can consult Lemma 5 of  \cite{as2} for a more self
contained proof of this fact, where the Clifford multiplication is
identified with the cross product operation.

\vspace{.1in}

\section { Associative Submanifolds }

\vspace{.1in}

 Any imbedding of a $3$-manifold   $f: Y^{3}\hookrightarrow M^7 $  induces
 an imbedding  $\tilde{f}: Y \hookrightarrow   \tilde{M}$:
\begin{equation}
\begin{array}{lcl}
  & & \tilde{M}\supset \tilde{M}_{\varphi } \\
  \hspace{.25in} \tilde{f}   \hspace{-.1in}&  \nearrow \;  &   \downarrow\\
\;Y\;\;& \stackrel{f}{\longrightarrow}   & M
\end{array}
\end{equation}

 \noindent and the pull-backs $\tilde{f}^{*} \Xi=T (Y)$ and  $\tilde{f}^{*} \V=\nu(Y)$ give the tangent and normal bundles of $Y$.  Furthermore,   if $f$ is an imbedding of an  associative submanifold into  a $G_2$ manifold  $(M, \varphi)$,  then the image of $\tilde{f}$ lands in $ \tilde {M}_{\varphi }$.  We will denote this canonical  lifting of any $3$-manifold $Y\subset M$ by $\tilde{Y}\subset \tilde{M}$.  Also since we have the dependency  $\V=\V_{\varphi }$, we can denote $\nu(Y)=\nu(Y)_{\varphi} =\nu_{\varphi}$ when needed.

\vspace{.1in}

$ \tilde {M}_{\varphi }$ can be thought of as a universal space
parameterizing associative submanifolds of $M$. In particular,  if
$\tilde{f}: Y\hookrightarrow  \tilde{M}_{\varphi}$ is the lifting
of an associative submanifold, by pulling back we see that  the
principal  $SO(4)$   bundle $\CP(\V)\to \tilde{M}_{\varphi} $
induces an $SO(4)$-bundle $\CP (Y)\to Y$, and gives the following
vector bundles  via the  representations:

\begin{equation}
\begin{array}{lcc}
\; \nu(Y) \hspace{.2in}:   &y \mapsto  q y\lambda^{-1}   & \hspace{1in}  \\
 \; T (Y)  \hspace{.1in}:  & x \mapsto  qx\ q^{-1} &   \\
%\;\lambda_{-}(\nu): \hspace{.4in}: & y  \mapsto  \lambda y  \lambda ^{-1} & \\
 \end{array}
 \end{equation}

where $[q,\lambda]\in SO(4)$, $\nu=\nu(Y)$ and $T(Y)= \lambda_{+}(\nu ) $. Also we can identify $T^{*}Y$ with $TY$ by the induced metric. From above we  have the action
$T^{*}Y \otimes \nu \to \nu  $   inducing  actions $\Lambda^{*}(T^{*}Y) \otimes \nu \to \nu  $.

\vspace{.1in}

Let $\L=\Lambda^{3}(\Xi)\to \tilde{M}$  be the determinant (real) line  bundle. Recall that the  definition (9) implies that  $\chi $ maps every oriented $3$-plane in $T_{x}(M)$ to its complementary subspace,  so $\chi$ gives a bundle map $\L \to \V$ over $\tilde{M}$, which is a section of $\L^{*}\otimes \V \to \tilde{M}.$  Since $ \Xi $ is oriented  $\L$ is trivial,   so $ \chi $ actually gives  a section

\begin{equation}
 \chi =\chi_{\varphi} \in \Omega^{0}(\tilde{M},  \V)
 \end{equation}

Clearly $ \tilde {M}_{\varphi }\subset \tilde{M}$  is the codimension $4$  submanifold which is the zeros of this section.  Associative submanifolds $Y\subset M$ are characterized by the condition $\chi|_{\tilde{Y}}=0$,  where $\tilde{Y}\subset \tilde{M}$ is the canonical lifting of $Y$. Similarly $\varphi $ defines a map $\varphi: \tilde{M}\to \R$.

\vspace{.1in}

\subsection { Pseudo-associative submanifolds}

$\;$

\vspace{.1in}

Here we generalize associative submanifolds to a more flexible class of submanifolds. To do this we first generalize the notion of imbedded submanifolds.

%\vspace{.1in}

{\Def A {\it Grassmann-framed  $3$-manifold}  in $(M,\varphi )$ is  a triple $(Y^{3} ,f, F)$, where $f:Y \hookrightarrow M$ is an imbedding, $F:Y\to \tilde{M} $, such that the following commute
 \begin{equation}
\begin{array}{rrc}
  & & \tilde{M}  \\
 F  \hspace{-.1in} &  \nearrow \;  &   \downarrow\\
\;\; Y\;\; & \stackrel{f}{\longrightarrow}   & M
\end{array}
\end{equation}
We call $(Y, f, F)$  a  {\it pseudo-associative submanifold } if
in addition $Image(F) \subset \tilde{M}_{\varphi}$. So a
pseudo-associative submanifold $(Y, f, F)$  with $F=\tilde{f}$ is
associative}.

{\Rm The bundle $ \tilde{M} \to M$ always admits a section, in
fact the subbundle $ \tilde{M}_{\varphi}  \to M$ has a section.
This is because  by \cite{t} every orientable $7$-manifold admits
a non-vanishing linearly independent $2$-frame field $\Lambda
=\{v_1, v_2 \}$\footnote {We thank T.Onder for pointing out
\cite{t}}. By Grahm-Schmidt process with metric $g_{\varphi }$, we
can assume that $\Lambda$ is orthonormal. The cross product
assigns $\Lambda$ to an orthonormal $3$-frame field $\{v_1,v_2,v_1
\times_{\varphi}v_2\}$ on $M$,  then  $ 3$-plane generated by $\{
v_{1}, v_{2},  v_1\times_{\varphi }v_{2}\}:= <v_{1}, v_{2}, v_{1}
\times_{\varphi }v_{2}>$}
 gives a section of  $ \lambda_{\varphi}: M\to \tilde{M}_{\varphi} $.
\begin{figure}[ht]  \begin{center}
\includegraphics{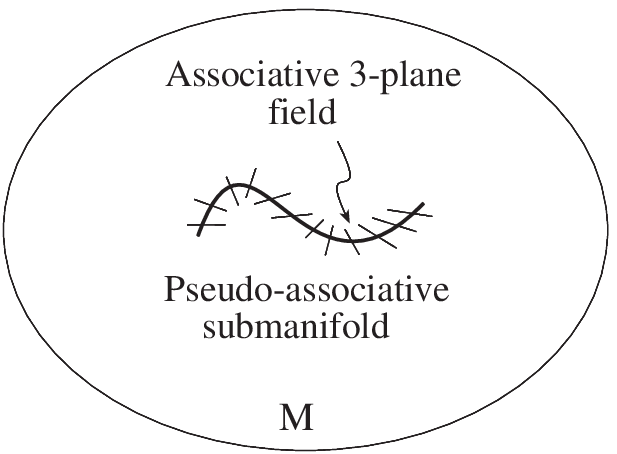}   \caption{}    \end{center}
\end{figure}
Let  $\CZ (M)$ and  $\CZ_{\varphi } (M)$ denote the set of Grassmann-framed   and the pseudo-associative submanifolds, respectively, and  let  $\CA_{\varphi }(M)$ be the set of associative submanifolds. We have inclusions
  $\CA _{\varphi} (M) \hookrightarrow \CZ_{\varphi }(M)\hookrightarrow \CZ(M) $,
  where the first map is given by  $(Y,f)\mapsto (Y, f, \tilde{f})$.
So  there is an inclusion  $Im (Y, M)\hookrightarrow \CZ(M)$, where  $Im(Y,M)$ is  the space of imbeddings.  This inclusion can be thought of  the canonical sections of a bundle
   \begin{equation}
\CZ(Y)   \stackrel{\pi}{\longrightarrow}  Im(Y, M)
\end{equation}
with fibers $\pi ^{-1}(f) = \Omega^{0}(Y, f^{*}\tilde {M})$. We
also have the subbundle $\CZ_{\varphi}(Y)
\stackrel{\pi}{\longrightarrow} Im(Y, M)$  with fibers
$\pi^{-1}(f) = \Omega^{0}(Y, f^{*}\tilde{M}_{\varphi})$.  So
$\CZ(Y) $ is  the set of  triples $(Y, f, F)$ (in short just set
of $F$'s),  where $F:Y\to \tilde{M}$ is a lifting of the imbedding
$f:Y\hookrightarrow M$. Also  $\CZ_{\varphi}(Y) \subset \CZ(Y) $
is a smooth submanifold, since $\tilde{M}_{\varphi}\subset
\tilde{M}$ is smooth. There is  the canonical section $\Phi: Im
(Y,M) \to \CZ(Y) $ given by $\Phi (f) = \tilde{f}$. Therefore,
$\Phi^{-1} \CZ_{\varphi}(Y):= Im_{\varphi}(Y, M)$  is the set of
associative imbeddings  $Y\subset M$. Also, any $2$-frame field
$\Lambda$ as above gives to a section $\Phi_{\Lambda}(f)=
\lambda_{\varphi}\circ f $. To make these definitions  parameter
free we also have to divide $Im (Y,M)$ by the diffeomorphism group
of $Y$. There are also  the vertical tangent bundles of  $\CZ(Y)$
and $\CZ_{\varphi} (Y)$
$$
\begin{array}{ccc}
T^{v}\CZ(Y) & \stackrel{\pi}{\longrightarrow} & \CZ(Y)\\
\cup &  & \cup\\
T^{v}\CZ_{\varphi} (Y) & \stackrel{\pi |}{\longrightarrow} & \CZ_{\varphi}(Y)
\end{array}
$$
with fibers
$\pi^{-1}(F)=\Omega^{0}(Y, F^{*}(\Xi^{*}\otimes \V))$.
 By Lemma 3 the fibers of $T^{v}(\CZ_{\varphi })$ can be identified with  the  kernel of the map induced by the Clifford multiplication
%\vspace{.02in}
\begin{equation}
 c: \Omega^{0}(Y,  F^{*}(\Xi^{*}\otimes \V))\to \Omega^{0}(Y, F^{*}( \V))
 \end{equation}
%\vspace{.01in}

 One of the nice properties of  a pseudo-associative submanifold $(Y,f, F) $ is that there is a Clifford multiplication action (by pull back)
  \begin{equation}
  F^{*}(\Xi^{*})\otimes F^{*}( \V) \to  F^{*}( \V)
  \end{equation}

\noindent If  $F$ is close to $ \tilde{f}$, by parallel
translating the fibers over $F(x)$ and $\tilde{f}(x)$  along
geodesics in $\tilde{M}$  we get canonical  identifications:
\begin{equation}
F^{*}(\Xi )\cong TY \;\;  F^{*}(\V )\cong \nu_{f}
\end{equation}
inducing Clifford multiplication between the tangent and the
normal bundles. So if $ \forall x \in Y$ the distance between
$F(x)$ and $\tilde{f}(x)$ is less then the injectivity radius
$j(\tilde{M})$, there is a Clifford multiplication between
the tangent and normal bundles of $Y$.

 \vspace{.1in}

 \subsection{ Dirac operator }

 $\:$

 \vspace{.1in}

The normal bundle $\nu=\nu(Y)$  of any orientable $3$-manifold $Y$
in a $G_2$ manifold $ (M,\varphi) $ has a $Spin(4)$ structure
(e.g. \cite{b2}). Hence we have $SU(2)$ bundles $S$ and $E$   over
$Y$  such that $\nu =S\otimes_{\H} E$ (18),  with $SO(3)$
reductions $ad S=\lambda_{+} (\nu)$, and  $ad E= \lambda_{-}(\nu)$
which is also  the bundle of endomorphisms $End (E)$.  If  $Y$ is
associative, then  the bundle $ad(S)$ becomes isomorphic to $TY$,
i.e. $S$ becomes  the spinor bundle of $Y$,  so $\nu (Y)$ becomes
a twisted spinor bundle.

\vspace{.1in}

The Levi-Civita connection of the $G_{2}$ metric of $(M,\varphi) $
induces  connections on  the associated bundles $\V$ and $\Xi$ on
$\tilde{M}$. In particular, it induces connections on the tangent
and normal bundles of any submanifold $Y^{3}\subset M$. We will
call  these connections  the {\it background connections}. Let
$\A_{0}$ be the induced connection on  the normal bundle
$\nu=S\otimes E$. From the Lie algebra decomposition  $
so(4)=so(3)\oplus so(3) $,  we can write $ \A_{0}= B_{0}\oplus
A_{0}$, where $B_{0}$ and  $A_{0}$ are connections on $S$ and $E$,
respectively.

\vspace{.05in}

Let   $\CA(E)$  and $\CA(S)$ be the  set of connections on  the
bundles $E$  and $S$. Hence  $ A\in \CA(E) $,  $ B\in \CA(S)$  are
in the form $A=A_{0} +a$, $B=B_{0}+b$,   where $a\in \Omega^{1}(Y,
ad \;E) $ and $b\in  \Omega^{1}(Y, ad\; S) $.  So $\Omega^{1}(Y,
\lambda_{\pm}(\nu))$ parametrizes connections on $S$ and $E$, and
the connections on $\nu $ are in the form $\A= B\oplus A$. To
emphasize the dependency on $b$ and $a$ we sometimes denote
$\A=\A(b,a)$, and $\A_{0}=\A(0,0)=A_{0}$.

 \vspace{.1in}

Now, let $Y^3\subset M$ be any smooth  manifold.  We can ex press
the covariant derivative $\nabla_{\A}: \Omega^{0}(Y, \nu)\to
\Omega^{1}(Y, \nu)$ on $\nu$  by $ \nabla_{A}=\sum e^i \;\otimes
\nabla_{e_i}$, where $\{e_{i}\}$ and $\{e^{i}\}$ are  orthonormal
tangent and cotangent  frame fields of $Y$,  respectively.
Furthermore, if $Y$ is an associative submanifold,  we can use the
Clifford multiplication of (25) (i.e. the cross product)  to form
the twisted Dirac operator  $\Di _{\A}: \Omega^{0}(Y,\nu) \to
\Omega^{0}(Y,\nu)$
  \begin{equation}
  \Di _{\A}=\sum e^i \;. \nabla_{e_i}
 \end{equation}
The sections lying in  the kernel of this operator are usually
called  harmonic spinors twisted by $(E, \A)$. Elements of  the
kernel of $\Di_{A_{0}}$ are called the  harmonic spinors twisted
by $E$,  or just the twisted harmonic spinors.

 \vspace{.1in}

 \section{ Deformations }

 \vspace{.1in}

In \cite{m}, McLean  showed that  the space of associative
submanifolds of a $G_2$ manifold $(M,\varphi)$, in a neighborhood
of a fixed associative submanifold $Y$,  can be identified with
the  harmonic spinors on $Y$ twisted by $E$. Since the  cokernel
of the Dirac operator can vary,  the dimension of its kernel  is
not determined (it has zero index since $Y$ is  odd dimensional).
We will remedy this problem by deforming $Y$  in a larger class of
submanifolds. To motivate our aproach we will first  sketch a
proof of  McLean's theorem (adapting  the  explanation in
\cite{b3}). Let  $Y\subset M$ be an associative submanifold, $Y$
will determine a lifting $\tilde{Y} \subset \tilde{M}_{\varphi }$.
Let us recall that the  $G_2$ structure $\varphi $ gives a metric
connection on $M$, hence it gives a connection $A_{0}$ and a
covariant differentiation  in  the normal bundle $\nu(Y)=\nu$
$$\nabla_{A_{0}}: \Omega^{0}(Y, \nu)\to \Omega^{1}(Y, \nu)=\Omega^{0}(Y, \;T^{*}Y\otimes \nu)$$
Recall that we  identified   $T^{*}_{y }(Y) \otimes \nu_{y}(Y)$ by the tangent space of the Grassmannian of $3$-planes $TG(3,7)$ in $T_{y}(M)$. So the  covariant derivative lifts normal vector fields $v$  of $Y\subset M$  to vertical vector fields $\tilde{v}$  in $T(\tilde{M})|_{\tilde{Y}}$. We want the normal vector fields $v$ of $Y$ to move $Y$  in the class of associative submanifolds of  $M$,  i.e. we want the liftings $\tilde{Y}_{v}$ of the nearby copies $Y_{v}$  of $Y$  (pushed off by the vector field $v$) to  lie in
$\tilde{M}_{\varphi }\subset \tilde{M} $ upstairs, i.e. we want the component of $\tilde{v}$  in the direction  of  the normal bundle $\tilde{M}_{\varphi}\subset \tilde{M}$  to vanish. By Lemma 3,  this means $\nabla_{A_{0}}(v)$ should  be in the kernel of the Clifford multiplication $ \; c=c_{\varphi}: \Omega^{0}(T^{*}(Y)\otimes \nu)\to \Omega^{0}(\nu)$, i.e.
$\Di_{A_{0}}(v)=c (\nabla_{A_{0}}(v))=0$, where $\Di_{A_{0}}$ is the Dirac operator induced by the  background connection $A_{0}$, i.e. the composition
\begin{equation}
\Omega^{0}(Y, \nu) \stackrel{\nabla_{A_{0}}}{\longrightarrow } \Omega^{0}(Y, T^{*}Y\otimes
\nu)\stackrel{c}{\to}  \Omega^{0}(Y, \nu)
\end{equation}

The condition $\Di_{A_{0}}(v) =0$ implies $\varphi $ must be integrable at $Y$, i.e.  the $so(7)$-metric connection $\nabla_{A_{0}}$ on $Y$  coincides with  $G_2$-connection (c.f. \cite{b2}).

\vspace{.05in}

Now we give a general version  of   the McLean's theorem,  without
integrability assumption on  $\varphi$: Recall from (Section 3.1)
that $\Phi^{-1} \CZ_{\varphi}(Y)$  is  the set of associative
submanifolds $Y\subset M$, where $\Phi: Im (Y,M) \to \CZ(Y) $  is
the canonical section (Gauss map) given by  $\Phi (f) =
\tilde{f}$. Therefore,  if  $f :Y \hookrightarrow M$ is  the above
inclusion, then  $\Phi(f) \in \CZ_{\varphi }$. So this moduli
space is smooth if  $\Phi $ was transversal to $\CZ_{\varphi}(Y)$.

\begin{figure}[ht]  \begin{center}
\includegraphics{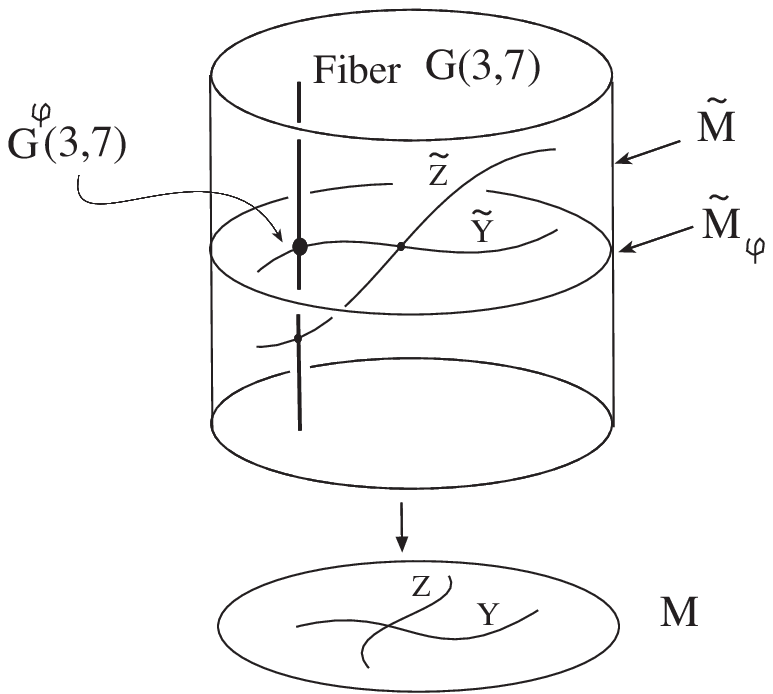}   \caption{}    \end{center}
\end{figure}

{\Thm  Let $(M^7, \varphi)$ be a manifold with a $G_2$ structure,
and $Y^3\subset M$ be an associative submanifold. Then the tangent
space of associative submanifolds  of $M$ at $Y$ can be identified
with the kernel of  a Dirac operator $\Di_{A}:\Omega^{0}(Y,
\nu)\to \Omega^{0}(Y, \nu)$, where $A=A_{0} +a$, and $A_0 $ is the
connection on $\nu$ induced by the metric $g_{\varphi}$,   and
$a\in \Omega^{1}(Y, ad (\nu )) $. In the case  $\varphi$ is
integrable $a=0$. In particular, the space of associative
submanifolds  of $M$ is smooth at $Y$ if the cokernel of $\Di_{A}$
is zero. }

\proof Let $f: Y\hookrightarrow M $ denote the imbedding. We
consider unparameterized deformations of  $Y$ in $ Im(Y,M) $ along
its normal directions. Fix a trivialization $ TY\cong im (\H) $,
by (17) we have an identification $ \tilde{f}^{*}(T^{v}\tilde{M})
\cong TY^{*}\otimes \nu + TY +  \nu $. We first claim  $\Pi \circ
\; d\Phi (v)=\nabla_{A}(v) $, where $d\Phi$ is the induced map on
the tangent space and
 $\Pi$ is the vertical projection.
\[
\begin{array}{rcl}
 \Omega^{0}(Y, \nu)=T_{f} Im (Y,M)  &  \stackrel {d\Phi } {\longrightarrow }     &
T_{\tilde{f}}\CZ(Y)= \Omega^{0}(Y, \tilde{f}^{*}(T^{v}\tilde{M}) ) \stackrel{\Pi}{\to}\Omega^{0}(Y, T^{*}Y\otimes \nu) \\
   &   &   \\
 \downarrow  exp   &   &   \downarrow  exp \\
  &   &   \\
  Im(Y,M) & \stackrel {\Phi } {\longrightarrow }    & \CZ(Y)
\end{array}
\]

\noindent The two vertical maps  $v\to f_{v}$, and $w\to  (\tilde{f})_{w}$ are  exponential projections of  tangent vectors,   i.e.    $f_{v}(y)=exp_{f (y) } (v)$ and   $(\tilde{f})_{w}(y)=exp_{\tilde{f}(y) } (w)$.
It suffices to check this claim  pointwise. Here for convenience view $f$ as an inclusion $Y\subset M$.

\vspace{.05in}
Let $ y=(y_1,y_2,y_3) $ be the normal coordinates of $Y$ centered around $y_{0}$, and  $\{e_{j}\}_{j=1}^{7}$ be an orthonormal frame field of $M$ defined on $Y$, with
$e_{j}(y_{0})=\partial/\partial y_{j} $ for $j=1,2,3$. To this data we can associate {\it Fermi coordinates}
$(y,t)$ around $f(y_{0})\in M$  (they are a version of normal coordinates along a submanifold, see for example \cite{g}):
\begin{equation}
(y,t )\longleftrightarrow f_{\sum t_{\alpha} e_{\alpha }  } (y)
\end{equation}
where  $t=(t_4,..,t_7)$.
Then  we can write
 $\tilde{f}(y_{0})= e_{1} \wedge e_{2} \wedge e_{3} $.  Hence by  definition we can express $d\Phi (v) =
 \widetilde{({f}_{v} )}=(f_{v})_{*}(e_{1})\wedge (f_{v})_{*}(e_{2})\wedge (f_{v})_{*}(e_{3}):=
 e_{1}(v)\wedge e_{2} (v) \wedge e_{3}(v) $.
 \begin{equation} d\Phi (v) (y_{0})=
\CL_{v} (e _1 \wedge e _2\wedge e_3)=\sum_{i=1}^{3} (*e _j) \wedge  \CL_{v}(e _j)|_{Y}\end{equation}

\noindent where $\CL_{v}$ denotes Lie derivative along $v$, and $*$ is the star of $Y$.
The metric connection is torsion free hence
$\CL_{v}(e _j)=\bar { \nabla}_{e_{j}}(v)-\bar{\nabla}_v(e_j)$, where $\bar{ \nabla}$ is the metric connection of $M$.  In case $(M,\varphi)$ is a $G_2$ manifold (i.e. when $\varphi$ integrable), by  (2) and (7), up to quadratic term  $\varphi$ is  $\varphi_{0}$, therefore we can write:
$$0=  {\bar \nabla}_{v}(\varphi )|_{Y} =\bar{ \nabla}_{v}(e^1\wedge e^2\wedge e^3)|_{Y}=
\sum_{j=1}^{3}(*e^{j})\wedge \bar{ \nabla}_{v}(e^{j}) |_{Y}\;,\;\; \mbox{which implies}$$
\begin{equation}
\Pi\circ d\Phi (v) (y_{0})= \sum_{j} (*e _j) \wedge  \nabla_{e_j}(v)
\end{equation}
where where $\{e^{j}\}$ is the dual coframe, and
$\nabla_{e_{j}}(v) $ is  the normal component of
$\bar{\nabla}_{e_{j}} (v)$, i.e. it is the induced connection on
$\nu(Y)$.  The expression (38) can be viewed as an infinitesimal
deformation of the $3$-plane $\tilde{f}(y_{0})$.  By the
identification $  *e _j \leftrightarrow e^{j} $ we can view it as
an element of the tangent space $T^{*}Y\otimes \nu $  of the
Grassmannian of $3$-planes in $T_{y_0}(M)$
 \begin{equation}
 \Pi\circ d\Phi (v) (y_{0}) =\sum e^{j}  \otimes \nabla_{e_{j}}(v)(y_{0})  =\nabla_{A_{0}}(v)(y_{0})
 \end{equation}

 \vspace{.1in}

 When $\varphi$ is not integrable, there is an extra term which we can write
$$\sum (*e_{j})\wedge \bar{\nabla}_{v}(e_{j})= \sum (*e_{j})\wedge \nabla_{v}(e_{j})$$
where $ \nabla_{v}(e_{j})$ is the normal component of $ \bar{\nabla}_{v}(e_{j})$. Notice $ <\bar{\nabla}_{v}(e_{k}),e_{k}>=0$, which is implied by  $v<e_k,e_k>=0$. So in this case (39) becomes
\begin{equation}
\Pi\circ d\Phi (v) (y_{0}) = \nabla_{A_{0}}(v) +a(v) = \nabla_{A}(v)
\end{equation}
where  $a(v) =
\sum e^{j}\otimes \nabla_{v}(e_{j})\in \Omega^{1}(Y, ad(\nu))$
and $A=A_{0}+a$. It easy to check that the expression $a(v)$ is independent of the choice the orthonormal frame $\{e_j\}$.

 \vspace{.1in}

By (31) the vertical tangent space of $\CZ_{\varphi} (Y)$ is given
by the kernel of the Clifford multiplication $c_{\varphi}:
\Omega^{0}(T^{*}Y\otimes \nu)\to \Omega^{0}(\nu)$. So,  locally
the moduli space of associative submanifolds of $(M,\varphi )$ is
given by the kernel of  $\Di_{A}$, i.e. the condition that $d\Phi
(v)$ lies in  $T^{v}_{f} \CZ_{\varphi} (Y)$  is given by
$\Di_{A}(v)=0$.  The moduli space is smooth if $\Phi $ is
transversal to $\CZ_{\varphi }(Y)$, i.e.
  if  the cokernel of $\Di_{A}$ is  zero. Since $T^{v}_{\tilde{f}} \CZ(Y) =\Omega^{0}(T^{*}Y\otimes \nu)$ and
 $$T\CZ(Y)/T\CZ_{\varphi}(Y)=  T^{v}\CZ(Y)/T^{v}\CZ_{\varphi}(Y)$$
  to check transversality  we look at the  induced maps, and use  $\Pi \circ  d\Phi_{f}(v)=\nabla_{A}(v)$
    \begin{equation*}
\Omega^{0}(\nu)=T_{f} Im (Y,M)  \stackrel {d\Phi } {\longrightarrow }  T_{\tilde{f}}\CZ(Y)
\stackrel {\Pi } {\longrightarrow } T^{v}_{\tilde{f}} \CZ(Y) \supset
T^{v}_{\tilde{f}} \CZ_{\varphi} (Y) \;\;\;\;\;
\mbox{\qed}
\end{equation*}

{\Rm  This theorem can also be proved by generalizing  McLean's
proof: The condition that an associative $Y\subset M$ remains
associative, when moved via the exponential map along a normal
vector field $v\in \Omega^{0}(Y, \nu ) $,  is $ \CL_{v}(\chi
)|_{Y}=0 $. We can choose local coordinates $(x_1,..., x_7)$ on
$M$,  such that $(x_1,x_2,x_3)$ gives the coordinates of $Y$. By
(13) and (14) $\chi= \sum  a^{\alpha}_{J}  \; dx^{J} \otimes
\partial /\partial x_{\alpha} $, with $\alpha \notin J$ and
$a^{\alpha }_{123 }|_{Y}=0$
%\frac{\partial}{\partial x_{\alpha}} $ and

\begin{equation}
\CL_{v}(\chi )|_{Y}=\sum  v(a^{\alpha}_{123}) \frac{\partial}{\partial x_{\alpha}}+
\sum  a^{\alpha }_{J } \; \CL_{v}(dx^{J}) |_{Y}
\otimes \frac{\partial}{\partial x_{\alpha}} =0
\end{equation}

McLean treated  integrable  $\varphi$ case, i.e  when
$(M,\varphi)$ is a $G_2$ manifold. In this case  the first term
vanishes, and  the second becomes $\Di_{A_{0}}(v)\otimes
dx^{123}$. But notice that the first term $a(v)$  is linear in $v$
and takes values in $\Omega^{0}(Y, \nu)$,  hence
$a\in\Omega^{0}(Y, ad (\nu) )$. So in the non-integrable case we
get a twisted  Dirac equation $ \Di_{A}(v)=0$, where $A=A_{0}+a$.
}

\vspace{.05in}

If in the proof of Theorem 4 we replace the $G_2$ structure
$\varphi $ with another $G_2$ structure  $\psi $ inducing the same
metric, the identification of the bundle $TY^{\perp}=\nu$ doesn't
change but the Clifford action $c_{\varphi}$ changes to another
one $c_{\psi}$,  corresponding to another $4$-dimensional Clifford
representions of  $T^{*}Y$.  These two representations  are
conjugate by a gauge automorphism $\gamma$ of $\nu$.
\begin{equation*}
\begin{array}{ccccc}
  & &\Omega^{0}(T^{*}Y\otimes \nu )&\stackrel{c_{\psi}}
   {\longrightarrow}& \Omega^{0}(\nu ) \\
   \\
&  &1\otimes \gamma \; \downarrow & &  \gamma \; \downarrow   \\
\\
 & &\Omega^{0}(T^{*}Y\otimes \nu ) & \stackrel{c_{\varphi}}
  {\longrightarrow} & \Omega^{0}( \nu )
 \end{array}
 \end{equation*}

\noindent Therefore,  if we call the Dirac operator induced by
$\psi $ by $\Di_{A_{1}}$, we  can write
$$\gamma (\Di_{A_{1}}(w)) =\sum dy^{j}.\gamma(\nabla_{j}(w)) = \sum dy^{j}.(\nabla_{j}\gamma(w))-dy^{j}.(\nabla_{j}\gamma)(w) $$
where the dot $''. ''$ denotes the Clifford product $c_{\varphi}$.
So, $D_{A_{1}}(w)=0$ gives a twisted version of the Dirac equation
$D_{A_{0}}(v)=0$ where $v=\gamma(w)$, this is because $\gamma
(\Di_{A_{1}}(w))=\Di_{A_{0}+a} (\gamma(w))$, where $a=-\sum dy^{j}
. (\nabla_{j}\gamma) \gamma^{-1}$. In Theorem 6 we will use the
twisting of the Dirac operator, under deformations of $\varphi$,
to obtain its surjectivity.

\vspace{.1in}

\section {Transversality }

\vspace{.1in}

We  can make the cokernel of Dirac operator $\Di_{A_{0}}$ zero  either by deforming the Gauss map $\Phi: Im (Y,M) \to \CZ(Y) $, or  by deforming the $G_2$ structure $\varphi$. Changing $\varphi $ can be realized by deforming $\varphi $ by a gauge transformation: Recall that the $G_2$ structures $\varphi $ on $M$  are the sections $\Omega^{3}_{+}(M)$ of the bundle (5). Also $GL(7,\R)$ conjugates $G_{2}=G_2^{\varphi_{0}}$ to any other $G_{2}$ subgroup
$G_{2}^{\varphi} $ of  $GL(7,\R)$ where
$$G_{2}^{\varphi}  =\{ A\in GL(7,\R) \;|\; A^{*}\varphi =\varphi\} \stackrel{\varphi}{\hookrightarrow}SO(7)$$

If we are interested in the $G_2$ structures inducing the same
metric, we  replace  $GL(7,\R)$ with $SO(7)$.  $SO(7)$ acts  on
$G(3,7)$  permuting submanifolds $G^{\varphi}(3,7)$,  where
$\varphi \in \Omega^{3}_{+}(M)$. More generally the gauge group
$\CG (P)$ of $P=P_{SO(7)}\to M$ acts on $\tilde{M}$ permuting
$\tilde{M}_{\varphi}$'s. Recall that $\CG (P)=\{
P\stackrel{s}{\to} P\;|\: s(pg)=s(p)g \;\}$, which can be
identified with sections  $\Omega^{0}(M; Ad(P))$ of the bundle
$Ad(P)\to M$ (c.f. \cite{amr}), where
$$  Ad(P)=P\times_{Ad} SO(7)=\{[p,h]\; |\;(p,h)\sim(pg,g^{-1}hg)\}$$
 $$\hbox{One can also  identify:}\;\;\;\CG (P)=\{ s:P \to SO(7) \;|\; s(pg)=g^{-1}s(p)g\;\}$$

\vspace{.1in}

\noindent The tangent space of  $\CG (P)$ at the identity  $I$ are the sections $ \mathfrak{ g} (P)= \Omega^{0}(M, ad(P))$  of the associated bundle  of  Lie algebras  $ ad(P)=P\times_{Ad} so(7) \to M$. Similarly $$  \mathfrak{ g} (P)=\{ h:P \to so(7) \;|\; h(pg)=h(p)g-gh(p)\;\}$$

We can  identify $T_{s}(\CG_{P}(M))\stackrel{\cong}{\mapsto}
 \mathfrak{ g} (P)$, by $s\mapsto s^{-1}ds$. There is also an  action $ \CG (P)\times \tilde{M} \to \tilde{M}$
given by $(s, [p,L])\to s[p,L]:= [ps(p),s(p)L]$, which we will simply denote it by
$ (s,L)\mapsto s.L $. There is the  pull-back action
 $\CG (P) \times \Omega^{3}_{+}(M) \to \Omega^{3}_{+}(M)$ given by
 $(s,\varphi)\to s^{*}(\varphi)$.  In particular, $s\tilde{M}_{\varphi}=\tilde{M}_{s^{*}\varphi}$.  Put another way,  if $\chi =\chi_{\varphi}$ is the $3$-form of Definition 3 and $L\in \tilde{M}$,  then  $\chi|_{L}=0 \iff s^{*}\chi |_{s^{-1}L}=0$. Hence   the $3$-plane $sL$ is $\varphi$-associative $\iff $ $L$ is $s^{*}\varphi$ -associative (similar to the process in the Kleiman transversality, c.f. \cite{ak})

 \vspace{.05in}

From the above action, we see that the space of $G_2$ structures
$\tilde{\Omega}^{3}_{+}(M)$ which induce the same metric on $M$
has the following identification:

\vspace{.05in}

{\Lem Let $\CG (P_{G_2} )$ be the stabilizer of the action of $\CG (P)$ on
$\tilde{\Omega}^{3}_{+}(M)$ (i.e.  the gauge transformations fixing $\varphi $) then: }
$$\tilde{\Omega}^{3}_{+}(M)=\CG (P)/\CG (P_{G_2})=\Omega^{0}(M, P\times_{SO(7)} {\R\P}^{7})$$
 \proof Clearly $\CG (P)$ acts transitively on $\tilde{\Omega}^{3}_{+}(M)$ with stabilizer $\CG (P_{G_2})$. To see the second equality, we identify the fibers of  the coset space with the fibers of $\tilde{\Lambda}^{3}_{+}(M)\to M$ by the map:
  $$\R\P^7=SO(7)/G^{\varphi}_{2}\to \tilde{ \Lambda}^{3}_{+}(\R^7)$$
  $G_{2}^{\varphi}s \mapsto  s^{*}\varphi $. The adjoint action of $SO(7)$ on $SO(7)$ moves cosets $$G_{2}^{\varphi}s \mapsto (g^{-1}G_{2}^{\varphi}g) g^{-1}sg=G^{g^{*}\varphi}g^{-1}sg$$ Hence by the obove identification,  on $\R\P^7$  it induces $\varphi \mapsto g^{*}\varphi $.
 \qed

\vspace{.1in}

Now we can deform the canonical section   $\Phi: Im(Y,M)\to \CZ(Y)$ by the map
 \begin{equation}
\tilde{\Phi }: \CG (P)\times Im(Y,M) \to \CZ(Y)
\end{equation}
  \noindent  $\tilde{\Phi}(s,f)=\Phi_{s}(f)=s\Phi(f)=s(\tilde{f})$,  that is  $\tilde{\Phi}(s,f)(y)=s(f(y)) \tilde{f}(y)$. Notice $ \CG (P) $ acts on the sections of the bundle $\CZ(Y)\to Im(Y,M)$.

   \vspace{.05in}

  { \Thm $\tilde{ \Phi } $  is transversal to $\CZ_{\varphi } (Y)$. Also $ \Phi_{s} $  is transversal to $\CZ_{\varphi } (Y)$ for a generic choice of  $s$, equivalently $\Phi $ is transversal to $\CZ_{s^{*}\varphi}(Y)$  for a generic $s$. }

\proof:  Let $\tilde{\Phi} (s,f)\in \CZ_{\varphi}(Y)$.  We can check transversality of $\tilde{\Phi} $ at $(s,f)$ by computing its derivative.  By the Leibnitz rule and Theorem 4 we can compute
\begin{equation*}
s^{-1}\circ \Pi\circ d \tilde{\Phi} (h,v): \mathfrak{ g}_{P}(M)\oplus \Omega^{0}(\nu ) \to T^{v}_{s(\tilde{f})}\CZ(Y) \to T^{v}_{\tilde{f}}\CZ(Y)=\Omega^{0}(T^{*}Y \otimes \nu)
\end{equation*}
where  $d \tilde{\Phi} \;(h,v)=s(f)\; [ \nabla_{A_{0}}(v) + s^{-1}ds(v)\tilde{ f} \;]$, and $v=f_{v}$ is the perturbation of the inclusion $f$.
 Observe that $ad(P)=End(TM)$, and the map $y\mapsto s^{-1}ds(v) \tilde{f}(y)$ is a vertical deformation of the $3$-plane $y\mapsto \tilde{f}(y)=T_{y}Y$, hence it is a section of the pull back of the vertical tangent bundle of $\tilde{M}\to M$ over $Y$,  i.e. an element
$a(v)\in T^{v}_{\tilde{f}}\CZ(Y)= \Omega^{0}(T^{*}Y \otimes \nu)$. More specifically,
if we decompose $s^{-1}ds(v)$ as an element of $so(7)$ on $T_{f(y)}(M)=T_{y}Y \oplus \nu_{y} (Y)$ in block matrices we can write:

\begin{equation} s^{-1}ds(v) \;  |_{ f(y)} =
\left(
\begin{array}{ccc}
  *&  -\alpha (v)^{t}    \\
 \alpha (v) & *      \\
\end{array}
\right)
\end{equation}

\noindent Because $\alpha(v)$ is  linear in $v$, we can view $\alpha \in \Omega^{1}(Y, ad \;\nu)$, therefore we can express $s^{-1} \Pi\circ d\Phi_{s}(v)=\nabla_{A_{0}}(v) + \alpha(v)=\nabla_{\A}(v)$ with $\A=A_{0}+\alpha$.  So  the transversality is measured by the cokernel of the twisted Dirac operator
$ c_{\varphi} (\nabla_{\A})=\Di_{\A}$, where $c_{\varphi}$ is the Clifford multiplication. Now by choosing $\alpha(v)$ we show that we can make $\Di_{\A}$ onto. This is because,  if $\Di_{A_{0}}$ is not already onto, we choose $0 \neq w\in im (\Di_{A_{0}})^{\perp}$. By self adjointness of the Dirac operator $\;  0=<\Di_{A_{0}}(v),w>= <v, \Di_{A_{0}}(w)>$, for all $v$. So $\Di_{A_{0}}(w)=0$, by analytic continuation  $w\neq 0 $ on an open set. Then
$w\in im (\Di_{\A})^{\perp}$ implies $ <c_{\varphi}(\alpha (v)), w>=0 $ and hence $w=0$, which is a contradiction. The last implication follows from by choosing $s$ in (43) we can get the full Lie algebra $so(7)$, and hence $v\mapsto a(v)$ is onto, and  the Clifford multiplication $c$ is onto (Lemma 3).

\vspace{.05in}

So we obtain a smooth manifold $\tilde{\Phi}^{-1}\CZ_{\varphi}(Y)$, and by choosing a regular value $s$ of the projection $\tilde{\Phi}^{-1}\CZ_{\varphi}(Y) \to  \CG_{P}(M)$ we get $\tilde{\Phi}_{s}^{-1}\CZ_{\varphi}(Y)$ smooth (note that the derivative of the projection is Fredholm). Clearly the condition that
$\Phi_{s}$ transversal to $\CZ_{\varphi}(Y)$ is equivalent to $\Phi$ being  transversal to  $\CZ_{s^{*}\varphi}(Y)$. \qed

\vspace{.1in}

Theorem 6 says that  the space of  $ s^{*}\varphi $ associative deformations of   an $\varphi$ associative submanifold $Y\subset M$, where $s\in  \CG_{P}(M)$, is a  smooth (infinite dimensional) manifold.
Infinitesimally these deformations correspond to the kernel of the twisted Dirac operator, twisted by the connections in the normal bundle $\nu (Y)$.  Define
\begin{equation}
\sigma:  \CG (P) \to \Omega^{0}(\tilde{M}, \Xi^{*}\otimes  \lambda (\V) )
\end{equation}

\noindent by $\sigma(s)(L)(v)=\alpha_{s}(v,L)\in \Xi^{*}\otimes
\V$, where $\alpha_{s}(v,L)$ is obtained by decomposing $s^{-1} ds
(v) \in \Omega^{0}(M, ad(P))$ on $TM=L\oplus L^{\perp}$ as an
element of $so(7)$

\begin{equation} s^{-1}ds(v) \;  |_{ L} =
\left(
\begin{array}{ccc}
  *&  -\alpha(v,L)^{t}    \\
 \alpha(v, L) & *      \\
\end{array}
\right)
\end{equation}

\noindent We can think of $\Omega^{0}(\tilde{M}, \Xi^{*}\otimes
\lambda (\V)) $  as an universal space  parameterizing connections
on $\nu \to Y$. The Gauss map $\tilde{f}$ of any imbedding $f:
Y\hookrightarrow M$ pulls back $\Xi^{*}\otimes  \lambda (\V) $  to
the parameter space  $\Omega^{1}(Y, \lambda (\nu))$  of the
connections on $\nu(Y)$.
 \begin{equation}
 \Omega^{0}(\tilde{M}, \Xi^{*}\otimes \lambda (\V) ) \stackrel {\tilde{f}^{*}}{\longrightarrow}  \Omega^{1}(Y, \lambda (\nu) )
\end{equation}

\noindent Clearly the set $\Omega^{1}(\tilde{M},\lambda (\V))$ can
also be used as the universal parameter space. As in Section 3.2,
given any imbedding $f:Y\hookrightarrow M$,
 we can deform the background connection  $ A_{0} \to A= \A(b,a) $  in the normal bundle $\nu(Y)$,  with $b\in \Omega^{1}(Y,\lambda_{+}(\nu ))$ and $a\in \Omega^{1}(Y,\lambda_{-}(\nu ))$,  and get a perturbed version of (35)
\begin{equation}
\Omega^{0}(\nu) \times  \Omega^{1}(\lambda_{\pm}(\nu ))  \stackrel{\Di_{A}}{\longrightarrow} \Omega^{0}(\nu)
 \end{equation}

\noindent with the twisted Dirac equation  $\Di_{\A}(v)= c(\nabla _{\A}(v))=\Di_{A_{0}}(v) + \alpha v=0 $, where $\alpha=(b,a)$.
Here we prefer perturbing by $a$ (perturbing $b$ has the effect of perturbing the metric on  $Y$). A generic nonzero $a $ makes the  map
$v \mapsto \Di_{A_{0}+a}(v)$ surjective.  We can choose this perturbation term $a$ universally.

\vspace{.1in}

\section { Complex Associative Submanifolds }

 \vspace{.1in}

 Let $(M,\varphi)$ be a manifold with a  $G_2$ structure. Here we will study an interesting class of associative submanifolds whose normal bundles come with an almost complex structure. The  subgroups $U(2)\subset SO(4)\subset G_{2}=G^{\varphi}_{2}$, more specifically
$$ (S^{1}\times SU(2))/ \Z_{2}\subset (SU(2)\times SU(2))/\Z_{2}\subset G_{2}$$

\noindent  give a $U(2)$-principal  bundle
 $\CP_{G_2} (M)\to  \bar{M}_{\varphi }=\CP_{G_2}(M) /U(2) $. Also $ \bar{M}_{\varphi }$ is   the total space of  an $S^{2}$ bundle  $ \bar{M}_{\varphi}\to \tilde{M}_{\varphi }=\CP_{G_2}(M)/SO(4)$, which is just  the  sphere bundle
 \begin{equation}
 \bar{M}_{\varphi}=S(\Xi) \to \tilde{M}_{\varphi}
 \end{equation}
of the $\R^{3}$-bundle: $\lambda_{+}(\V)=\Xi \to \tilde{M}_{\varphi}$.
 We can identify the  sections  of  (48) with almost complex structures on $\V$. Notice $\bar{M}_{\varphi}\to M$ is a bundle with fibers
 $G_{2}/U(2)$, which we can view as the complex version of the associative Grassmanns $G^{\varphi}_{\C}(3,7)$.
In fact if $V_2(M)$ and $G_2(M)$ are the bundle of orthonormal $2$-frames and oriented $2$-planes in $M$ respectively,  the fibration  $V_2(M)\to G_2(M)$ can be identified by:
 \vspace{.03in}
     \begin{equation*}
 \CP_{G_{2}}(M)/SU(2)\to \CP_{G_{2}}(M)/U(2)
 \end{equation*}

  \vspace{.03in}

\noindent with its projection $\{u,v\}\mapsto ( <u,v, u\times v> , u\times v )$, and also the projection map $  \CP_{G_{2}} (M)\to \CP_{G_{2}}(M)/SU(2)$ on the fibers is given by the map $G_2\to V_{2}(\R^7)$ defined by $\{ v_1,v_2, v_3\}\mapsto \{v_1,v_2\}$ (recall the definition of $G_2$ in (1)).  Put another way, $G_2$ acts transitively on $V_{2}(\R^7)$  with stabilizer $SU(2)$. By summing up above:

{\Prop $\bar{M}_{\varphi}=S(\Xi)= G_{2}(M)$}

 \vspace{.1in}

 More generally, for the Riemannian manifold $(M^7, g_{\varphi})$ we can  take the sphere bundle $\bar{M}\to \tilde{M} $ of  $\lambda_{+}(\V)\to \tilde{M}$, and get codimension $4$ inclusion of the smooth manifolds $\bar{M}_{\varphi} \subset   \bar{M}$  (of dimensions $17 $ and $21$). The sections  of the bundle $\bar{M}\to \tilde{M}$ gives the parametrization of the almost complex structures on $\V$, and
 $\bar{M}\to M$ is a bundle with fibers $G_{\C}(3,7):=SO(7)/U(2)\times SO(3)$. For all $G_2$ structures $\varphi$ inducing the same metric on $M$, we have the inclusions $G_{2}^{\varphi}\hookrightarrow SO(7)$ inducing imbeddings $G_{2}(M)= \bar{M}{_{\varphi} } \hookrightarrow \bar{M}$, which is fiberwise  $<u,v>\mapsto <u,v,u\times v>$
$$G(2,7)= G_{2}^{\varphi}/U(2) \hookrightarrow SO(7)/U(2)\times SO(3)$$

\vspace{.1in}

By \cite{t} the bundle $V_{2}(M)\to M$ has always a section $\Lambda =\{u,v\}$, which induces sections of the bundles $\bar{M}_{\varphi}\to \tilde{M}_{\varphi}$ and $\tilde{M}_{\varphi}\to M$ (for simplicity we will abuse notation and denote all these sections by $\Lambda$ also). So $\Lambda $ gives an almost complex structure on $\V\to \tilde{M}_{\varphi}$. By $\Lambda$, we can pull back $\Xi$ and $\V$ to bundles ${\bf E}$ and ${\bf V}$ on $M$, respectively, and ${\bf V}$ has an almost complex structure (by the discussion following Lemma 2 we can describe this complex structure  with the cross product with $u\times v$).

\vspace{.1in}

{\Def From now on, we will denote a $7$-manifold with a $G_2$
structure and a nonvanishing $2$-frame field $\Lambda$ with
$(M,\varphi, \Lambda)$.}

\vspace{.1in}

Given $(M,\varphi, \Lambda)$, then  the induced $U(2)$ structure on $\V\to \tilde{M}_{\varphi}$ canonically lifts to a $Spin^{c}(4)$ structure  by the diagram:
\begin{equation}
\begin{array}{ccc}
  &   &   Spin^{c}(4)\\
  &  \nearrow &  \downarrow \\
 U(2) &  \to &   SO(4)\times S^1
\end{array}
\end{equation}

\noindent where
$ U(2)=(S^1\times S^3) /\Z_{2}$, $SO(4)=(S^3\times S^3)/\Z_{2}$, $Spin^{c}(4)=
(S^3\times S^3\times S^1)/\Z _{2}$,
where the horizontal map $[\lambda,A]\mapsto ([\lambda,A], \lambda^2)$ lifts to the map $[\lambda,A]\mapsto (\lambda,A, \lambda)$. This means there is a $\C^2$-bundle $\W\to \tilde{M}_{\varphi}$ with $\V_{\C}=\W\oplus \bar{\W}$, and transition function $\lambda^2$ gives the determinant line bundle $K=\Lambda^{2}\bar{\W} \to \tilde{M}_{\varphi}$. Also we can write $\V_{\C}=\W^{+}\otimes\W^{-}$ with $\W^{+}=K^{-1}+\C$ and $\W^{-}=\bar{\W}$. Recall  $\Xi^{*}=\lambda_{+}(\V)=\Lambda^{2}_{+}(\V)=K +\R$. We have a Clifford action $\V\otimes \W^{+}\to \W^{-}$ which extends  to  Cllifford action
$$\Xi^{*} \otimes \W^{+}\to \W^{+}$$
The identifications $\W^{+}\otimes \bar{\W}^{-}=\C\oplus \C\oplus K\oplus \bar{K}=(K\oplus \R)_{\C}\oplus \C =\Xi^{*}_{\C}\oplus \C$ gives
 the usual quadratic bundle map of the Seiberg-Witten theory (c.f. \cite{a}):
 \begin{equation}
 \sigma : \W^{+}\otimes  \W^{+} \to \Xi^{*} _{\C}
 \end{equation}
  \[
 \sigma(x,x)=  ( \frac{|z|^2 -|w|^2}{2},  \bar{z} w ) \;\;,\;\;\mbox{where}\;\;x=(z,w)\]

 % \vspace{.1in}

{\Def A submanifold $f: Y^3 \hookrightarrow M$  of  $
(M^7,\varphi, \Lambda)$  is called {\it $\Lambda$-associative} if
$\tilde{f}=\Lambda \circ f$ where $\tilde{f}$ is the Gauss map,
and it is called almost $\Lambda $-{\it associative} if it comes from
transverse section of the bundle ${\bf V}\to M$ (recall ${\bf V}$
is obtained from $\Lambda$).}

 \vspace{.1in}

\noindent An $\Lambda$-associative, more generally almost $\Lambda$-associative submanifold $Y$ of $(M,\varphi, \Lambda)$ induce canonical isomorphisms
 $TY\cong \tilde{f}^*(\Xi)$ and $\nu(Y)\cong \tilde{f}^*(\V)$ (by transversality). 
 \begin{equation*}
\begin{array}{rrc}
  & & \hspace{-.4in} \tilde {M}_{\varphi}  \\
 \tilde{f}  \hspace{-.1in}  &  \nearrow \;  &   \downarrow \uparrow \Lambda \;\; \;\;\;\;\;\;\; \\
\;\; Y\;\; & \stackrel{f}{\hookrightarrow}    &M \;\;\;\;\;\;\;\;\;
\end{array}
\end{equation*}

So normal bundle of any almost $\Lambda $-associative submanifold $Y^3 \subset M^7$ has a $U(2)$ structure, therefore it has a $Spin^{c}(4)$ structure, with induced $\C^2$-bundle $W\to Y$ and its determinant line bundle $K\to Y$, and  a Clifford action of $T^{*}Y\otimes W\to W $  (induced from the cross product). An example of a $\Lambda$ associative submanifold is the zero section of the spinor bundle $\CS\to Y^3$ (the $G_2$ manifold constructed in \cite{bsa}).

 \vspace{.05in}

In general, the background $SO(4)$ connection $A_{0}$ on the
normal bundle $\nu(Y)$ of a $\Lambda $-associative submanifold $Y$
may not reduce to a $U(2)$ connection if the $1$-form whose dual
gives the splitting $TY\cong K\oplus \R$ is not parallel.
Nevertheless from the $Spin^c(4)$ structure on $\nu(Y)$ we do get
a connection on the complex bundle $W\to Y$ provided we pick a
connection on the line bundle $K\to Y$ (from (49)). In the next
section we will study the local deformation space of $\Lambda
$-associative manifolds, by deforming them in the complex bundle
$W$, with the help of the connections on $K$.

%\vspace{.05in}

{\Rm  Associative submanifolds with $Spin^c(3)$ structure $(Y,c)\hookrightarrow (M, \varphi) $ come equipped with an  $U(2)$ structure (hence $Spin^c(4)$ structure) on their normal bundles (Lemma 2). We can free the deformations space these manifolds from the extra parameter $c$, by picking up  a generic $\Lambda$, and studying  the deformations of more relax almost $\Lambda $-associative submanifolds $Y\subset (M,\varphi, \Lambda)$. In this case the $Spin^c(3)$ structure on $TY$ comes from the pull-back. Also, by further deforming the $2$-plane field $\Lambda $ on $M$ deforms the %
the the $Spin^c$ structure on $Y$.}

{\Rm We could have considered complex structures on $\V$
corresponding to the right reduction, i.e. the subgroup  $
(SU(2)\times S^{1})/ \Z_{2}\subset (SU(2)\times
SU(2))/\Z_{2}\subset G_{2}$. In this case, they correspond to the
sections of the $S^2$ bundle
$\lambda_{-}(\V)\to\tilde{M}_{\varphi}$. Here we opted to the left
reductions since they concretely relate to $\Xi$ by
$\lambda_{+}(\V)=\Xi$.}

 \section{ Deforming $\Lambda $-associative submanifolds }

Let $(M,\varphi, \Lambda )$ be a manifold with $G_2$ structure and
a non-vanishing $2$-plane field, $\CM (M, \varphi ,\Lambda) $ be
the space of $\Lambda$-associative submanifolds. Here we will
study the local  ``complex'' deformations of  $\CM(M,\varphi
,\Lambda) $ near a particular $f: Y \hookrightarrow M$. These are
the deformations of $Y$  inside its complex normal bundle $W$,
with the help of the connections $\CA(K)$ on  the line bundle $K=
\mbox{ det }W$. These deformations are identified with the kernel
of a twisted Dirac operator twisted by  the connections in
$\CA(K)$. Introducing new variables  $\CA(K)$ makes the
deformation space smooth. Up to this point this section can be
viewed as a version of Theorem 4  for the $\Lambda$-associative
submanifolds. But now the connection parameter can be constraint
with the natural map (50) to obtain Seiberg-Witten like equations,
which gives a compactness result for this more restricted local
deformation space of $Y$. Reader should note that these equations
are $Spin^c(4)$ Seiberg-Witten equations on $Y^3$ (which are
usually associated to $4$-manifolds), as opposed to the usual
$Spin^c(3)$ Seiberg-Witten equations. The Clifford action
$T^*Y\otimes \W\to W$ is induced via the identification
$\Lambda_{+}^2W=T^*Y$, it is also induced by the cross product
operation on $M$.

\vspace{.05in}

Let $Y\in \CM (M, \varphi ,\Lambda) $.  Let $W\to Y$ be the
complex bundle associated to $\nu(Y)$, and $K\to Y$ be its
determinant  line bundle.  Let $B_{0 } $ be  the background
connection on $\nu(Y)$ (induced by $\varphi$), then as discussed
in last section $B_{0}$ along with  $A\in \CA(K)$ defines a
connection on $W\to Y$,  denote by $\A=   B_{0} \oplus A$. We can
write   $A=A_{0 }+a$ with $a\in \Omega^{1}(Y) =T_{A_{0}}\CA(K)$
(tangent space  of connections) and $\A= \A (a) $. Then we get a
complex version of the map (47) $\; (v,a)\mapsto \Di_{\A}(v)=
\Di_{\A(0)}(v) +a.v$
\begin{equation}
  \Omega^{0}(Y, W) \times \Omega^{1}(Y,i\R)  \stackrel{\Di_{\A}}{\longrightarrow} \Omega^{0}(Y, W)
 \end{equation}
which is  the derivative of a  similarly defined map
\begin{equation}
     \Omega^{0}(Y, W) \times \CA(K)   \to  \Omega^{0}(Y, W)
  \end{equation}
  
\begin{figure}[ht]  \begin{center}
\includegraphics{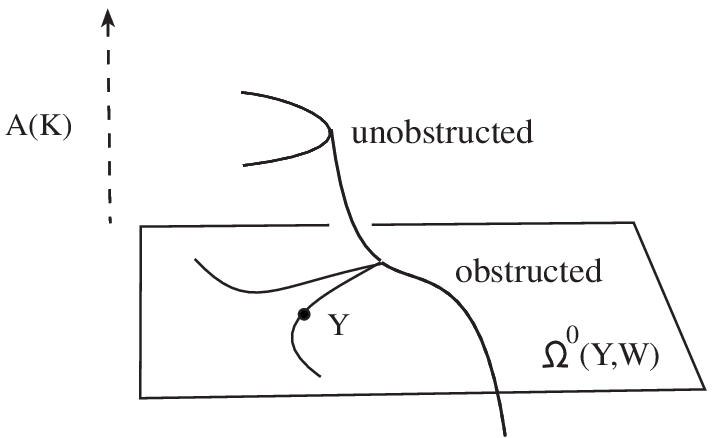}   \caption{}    \end{center}
\end{figure}

\vspace{.2in}

In each slice $\A(a)$,  we are deforming along normal vector fields by the connection $\A(a)$, which is a perturbation of the background connection $\A(0)$.    To get  compactness we can  cut down  this parametrized moduli space with an additional equation (induced from the map (50)) of the Seiberg-Witten theory  $\Psi^{-1}(0)$, where
\begin{equation*}
\Psi : \Omega^{0}(Y,W) \times \CA(K)\to \Omega^{0}(Y,W) \times \Omega^{2} (Y,i\R)
\end{equation*}
\begin{equation}
 \begin{array}{c}
 \Di_{\A}(v)=0 \\
*F_{A}=  \sigma(v,v)
 \end{array}
  \end{equation}

\noindent where $F_{A}$ is the curvature of the connection $A=A_0 +a$ in $K$, and $*$ is the star operator on $Y$. Note that $Y$ comes equipped with the natural submanifold metric.  Now we proceed exactly as in the Seiberg-Witten theory of $3$-manifolds (e.g. \cite{c}, \cite{lim}, \cite{ma}, \cite{w}). To obtain smoothness  of $\Psi^{-1}(0)$, we  perturb the equations by $1$-forms $\delta \in \Omega^{1}(Y)$ and get  a new  equation $\Phi=0$, where
\begin{equation*}
\Phi : \Omega^{0}(Y,W) \times \CA(K) \times  \Omega^{1} (Y) \to \Omega^{0}(Y,W)\times \Omega^{1} (Y,i\R)
\end{equation*}
\begin{equation*}
  \begin{array}{ccc}
 \Di_{\A}(v) &=& 0\\
 *F_{A} + i \delta &= & \sigma (v,v)
\end{array}
 \end{equation*}

\noindent We can choose the perturbation term universally $\delta ={f}^{*}(\Delta )$, where $\Delta \in \Omega^{1}(M )$.  Then  $ \Phi $ has a  linearization:
\begin{equation*}
D\Phi_{(v_0,A_{0},0)}: \Omega^{0}(Y,W) \times \Omega^{1}(Y,i\R) \times  \Omega^{1} (Y)
\to \Omega^{0}(Y,W)\times \Omega^{1} (Y,i\R)
\end{equation*}
$$D\Phi_{(v_0,A_{0},0)}\;(v,a,\delta)=(\Di_{A_{0}}(v) + a. v_{0}, \; *da + i\delta -2 \sigma(v_{0},v))$$

\noindent We see that  $\Phi^{-1}(0)$ is smooth and the projection $\Phi^{-1}(0) \to \Omega^1 (Y)$ is onto, so by Sard's theorem for  a generic choice of $\delta$ we can make $\Phi_{\delta}^{-1}(0)$  smooth, where  $\Phi_{\delta}(v,A)=\Phi (v,A,\delta)$. The bundle $W$ of $Y$ has a complex structure, so  the  gauge group  $\CG(K)=Map (Y, S^1)$ acts on the solution set  $\Phi^{-1}_{\delta}(0) $, and makes the quotient $\Phi^{-1}_{\delta }(0) / \CG(K)$ a  smooth zero-dimensional manifold. This is because  the infinitesimal   action of $\CG(K)$ on the complex
$\Phi_{\delta}: \Omega^{0}(Y,W)  \times  \CA(L)  \to \Omega^{0}(Y,W)\times \Omega^{1} (Y,i\R)$
 is given by the  map
$$\Omega^{0}(Y,i\R)\stackrel{G}{\longrightarrow} \Omega^{0}(Y,W)\times \Omega^{1}(Y,i\R)  $$
where $G( f )= ( fv_{0}, df)$.
So after dividing by $\CG$,  tangentially the complex $\Phi_{\delta}$ becomes
\begin{equation*}
\Omega^{0}(Y; i\R)\stackrel{G}{\longrightarrow}  \Omega^{0}(Y,W)\times \Omega^{1}(Y,i\R)
\to \Omega^{0}(Y,W)\times \Omega^{1} (Y,i\R)/G
\end{equation*}
 Hence  the index of this complex is the sum of the indices of the Dirac operator
 $\Di_{\A_{0}}: \Omega^{0}(Y,W)\to \Omega^{0}(Y,W)$ (which is zero), and the index of the following complex
 \begin{equation*}
 \Omega^{0}(Y,i\R)\times  \Omega^{1}(Y,i\R)
\to \Omega^{0}(Y,i\R)\times \Omega^{1} (Y;i\R)
\end{equation*}
given by $ (f,a)\mapsto (d^{*}(a), df+ *da) $, which is also zero  since $Y^3$ has zero Euler characteristic.
  Furthermore, $\Phi^{-1}_{\delta}(0) / \CG(K)$ is compact and oriented (the same proof as in the  Seiberg-Witten theory). Hence we get a number $SW_{Y}(M)$. Here we don't worry about metric dependence
of $SW_{Y}(M)$  since we have a fixed background metric induced
from the $G_2$ structure. Hence we associated a number to a
$\Lambda$-associative submanifold $Y$ of $(M,\varphi, \Lambda)$.
In particular, $Y$  moves in an unobstructed way along  the
parametrized  sections the complex normal bundle
$\Omega^{0}(Y,W)\times \CA(L)$. Furthermore all these constructions  work for almost $\Lambda$-associative submanifolds. So we have:

{\Thm Let $Y$   be an almost  $\Lambda$-associative submanifold of  $ (M,\varphi, \Lambda
)$. By cutting down the space of parametrized complex deformations
of $Y$  with an additional equation as in (53)  we obtain a zero
dimensional compact smooth oriented manifold, hence we can
associate a number $\Lambda_{\varphi} (Y)\in \Z$. }

{\Rm Clearly $\Lambda_{\varphi} (Y)$ is invariant under small isotopies
 through almost $\Lambda$-associative submanifolds $Y\subset (M,\varphi, \Lambda)$. }

\vspace{.1in}

The equations (53) can be induced  universally from equations on
$(M^7,\varphi, \Lambda)$  by restriction:  The $2$-frame field
$<u,v>$  gives a splitting of the tangent bundle  $TM=\bE\oplus
\bV$ with  an $SO(3)$ bundle $\bE=<u,v,u\times v>$  and a
$U(2)$-bundle   $\bV ={\bf E}^{\perp}$, such that $\lambda_{+}(\bf
V)=\bE$. Let $\bW\to M$ be the induced $\C^2$-bundle, and $\bK\to
M$ be the determinant line bundle of $\bW$.
  We can define an  action  $T^{*}(M)\otimes \bW \to \bW$:
  For  $ w=x+y \in TM $, with $ x\in \bE $, $ y\in\bV $ and $ z\in \bW $ with $w.z=xz$.
 It is easy to check that this is a partial Clifford action, i.e.
  $w.(w.z) =- |x|^{2} z$ Id, and it extends to an action $\Lambda^{2}(T^{*}M)\otimes \bW\to \bW$, and we have the map $\sigma: \bW\otimes \bW\to \bE_{\C}$  of (50).

 \vspace{.1in}

These bundles inherit connections from the Levi-Civita connection of $(M, g_{\varphi})$.
Let $\CA(\bK)$ be the connections on $\bK$. Let $A_{0}$ denote the background connections. Then any $A\in \CA(\bK)$ along with $A_{0}$ determines a connection  on $\bf W$.  Write $A=A_{0}+a$ with $a\in \Omega^{1}(M)$.
Hence for $A\in \CA(\bK)$ we can define  a partial Dirac  operator $\Di_{A}(v)=\Di_{A_{0}}(v) +a.v $ on $\bW\to M$,  which is the composition:
$$ \Omega^{0}(M,\bW)\stackrel{\nabla_{A}}{\longrightarrow } \Omega^{0}(M,T^{*} M \otimes \bW)
\stackrel{c_{\varphi}}{\longrightarrow }  \Omega^{0}(M,\bW)  $$

 We can now write the global version of the equations (54) on  $M$ in the usual way
 $$\phi: \Omega^{0}(M, \bW)\times  \CA(\bL)  \to   \Omega^{0}(M,  \bW)  \times  \Omega^{1} (M ) \;\;\;\mbox{which is}$$
 \begin{equation}
 \begin{array}{c}
 \Di_{A}(v)=0 \\
*F_{A} = \sigma(v,v)
 \end{array}
  \end{equation}
where $*:TM\to TM$ is the star operator on $\bE$ and zero on
$\bV$. We can perturb these equations by $1$-forms to $\Phi =0$,
and proceed as before.
 $\bW$ has a complex structure. The  gauge group  $\CG(\bL)=Map (M, S^1)$ acts on the solution set  $\Phi^{-1}(0) $,  and the quotient $\Phi^{-1}(0) / \CG(\bL )$ can be formed. To sum up  we have:

{\Prop  Any almost $\Lambda$-associative submanifold $f: Y^{3}\hookrightarrow (M,\varphi, \Lambda)$  pulls back the equations (54) to the Seiberg-Witten equations (53) on $Y$.}

\vspace{.1in}

\section{ Associative $3$-Plane Fields of $G_2$ Manifolds }

Recall that, any non-vanishing oriented $2$-plane field $\Lambda
=<u,v>$ on $(M,\varphi)$  determines a section $\Lambda_{\varphi}:
M\to \tilde{M}_{\varphi} \subset \tilde{M}$. In particular, it
gives  a non-vanishing associative $3$-plane field
$\bE=\bE_{\Lambda,\varphi}\to M$ on $M$, and a complex structure
on the complementary $4$-plane field
$\bV=\bV_{\Lambda,\varphi}\to M$, and  a splitting
$TM=\bE\oplus \bV $, with  $\lambda_{+}(\bV)=\bE$. From the
construction we get a further splitting $\bE={\bf \Lambda } \oplus
{\bf \xi}$, corresponding  to $<u,v>\oplus <u\times v>$. The
orientation of the $2$-dimensional bundle ${\bf \Lambda }$ gives
it a complex structure, and we have
 \begin{equation}
 TM=\bar{\bE}\oplus \xi
 \end{equation}
where $\bar{\bE}={\bf \Lambda} \oplus \bV $  is a $6$-plane bundle
with a complex structure and $\xi$ is the line bundle $<u\times
v>$. Note that if  $\varphi $ is integrable and the vector field
$u\times v$ is parallel then $M$ would be a Calabi-Yau $\times
S^1$ (since $G_2$ holonomy would reduce to $SU(3)$). So
non-vanishing oriented $2$-plane fields may be thought of  objects
taming the $G_2$ structure. Any integral submanifold of the
corresponding distribution $\bE$ is an associative submanifold
$Y^3 \subset M$ with a  $Spin^{c}$-structure (i.e. the $2$-plane
field $\xi=\Lambda|_{Y}$).

\begin{figure}[ht]  \begin{center}
\includegraphics{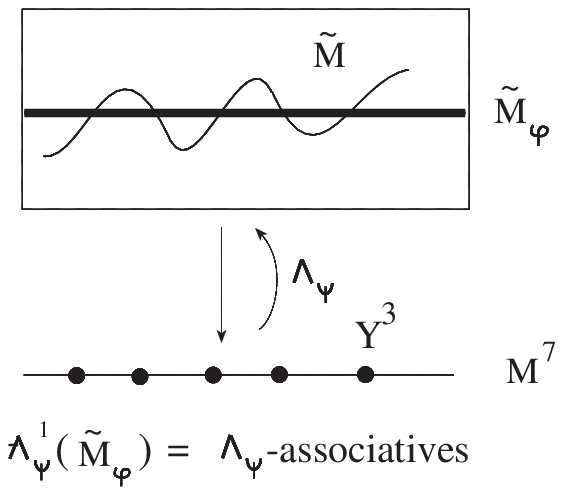}   \caption{}    \end{center}
\end{figure}

\vspace{.05in}

By fixing the plane field $\Lambda$, and varying $\varphi \in \tilde{\Omega}^{3}_{+}(M) $ (the set of $G_2$ structures inducing the same metric  on $M$) has the effect of  varying $\xi \in \Lambda^{\perp}$ (the cross product operation on $\Lambda$)  and varying the complex structure on $\bV=(\Lambda \oplus \xi )^{\perp}$. These $\xi$'s are the sections of the $S^4$-sphere bundle of  $\Lambda^{\perp}\to M^7$, hence generically any other section will agree with $\Lambda_{\varphi}$ on some $3$-manifold $Y\subset M$. We will show that this $3$-manifold is almost $\Lambda$-associative.  First consider the parametrized section:
\begin{equation}
{\bf \Lambda}: \tilde{\Omega}^{3}_{+}(M)\times M\to \tilde{M}
\end{equation}
 $(\lambda , x)\mapsto \Lambda_{\lambda }(x)$. By Lemma 5 there is an identification
$\tilde{\Omega}^{3}_{+}(M)=\{s^{*}(\varphi) \;|\; s\in \CG(P)\}$ (the sections of  an $\R\P^7$ bundle over $M$). We claim ${\bf \Lambda}$  is  transversal to $\tilde{M}_{\varphi}$.
 
\vspace{.05in}
First we need to recall a few facts: By \cite{b2}, the deformations of the $G_2$ structure $\varphi$  fixing the metric $g=g_{\varphi}$, are parametrized by $\varphi_{\lambda}$ below, where  $\lambda=[a,\alpha]$ are the sections  $\tilde{\Omega}^{3}_{+}(M)$ of  the $\R\P^7$-bundle, which is the projectivization $P(\R \oplus T^{*}M)\to M$ \begin{equation*}
\varphi_{\lambda}=(a^2-|\alpha|^2)\varphi +2a * (\alpha \wedge \varphi)+
2\alpha\wedge*(\alpha\wedge *\varphi)
\end{equation*}
where $a^2+|\alpha|^2=1$. From the identities $*(\alpha\wedge \varphi)=-\alpha^{\#} \lrcorner \; *\varphi $ and $*(\alpha\wedge *\varphi)=\alpha^{\#}\lrcorner \;\varphi $, where $\alpha^{\#}$ is the metric dual of $\alpha $,  we can also express
\begin{equation} \varphi_{\lambda}=\varphi -2 \alpha^{\#}\lrcorner \;[\;a(*\varphi)+\alpha\wedge \varphi \;]
\end{equation}
$$*\varphi_{\lambda}=*\varphi +2\alpha \wedge [\; a \varphi- (\alpha^{\#}\lrcorner \;*\varphi) \;]$$

\noindent  Not to clutter notations, we denote
$\Lambda_{\lambda}=\Lambda_{\varphi_{\lambda}}$ and use the metric
to identify $T^{*}(M)={\bf E}\oplus{\bf V}$, and identify $M$ with
the zero section of the bundle ${\bf V}\to M$.

{\Thm For $\alpha\in \Omega^{1}(M)$ which is a transverse section
of ${\bf V}\to M$, the map  ${\Lambda}_{\lambda }$, where $\lambda
=[a,\alpha ]$ and $a\neq 0$, is transversal to
$\tilde{M}_{\varphi}$, and
$\Lambda^{-1}_{\lambda}(\tilde{M}_{\varphi} )=\alpha^{-1}(M)$.

\proof

The set  ${\Lambda}_{\lambda}^{-1}(\tilde{M}_{\varphi})$  is given
by the solutions of the equation $\Lambda_{\varphi_{\lambda}}(x)=
\Lambda_{\varphi}(x)$, where  $\varphi \mapsto \varphi_{\lambda}$
is a deformation of $\varphi $.  Since ${\bf E}$ is obtained from
the oriented $2$-plane field $\Lambda =<u,v>$ by association
$<u,v>\mapsto <u,v,u\times_{\varphi}v>$, this equation is
equivalent to  $(u \times v)_{\lambda}(x)=(u\times v ) (x)$ (up to
positive scalar multiple), where $(u \times v)_{\lambda}$ denotes
the cross product corresponding to $\varphi_{\lambda}$. By using
$(u\times v)^{\#}=u\lrcorner \; v \lrcorner \;\varphi $, we can
calculate the deviation of the cross product operation under the
deformation
  \begin{eqnarray*}
(u \times v)_{\lambda}&= &(1-2|\alpha|^2)(u\times v) +
2\;{\bf [}-a \chi(u,v,\alpha^{\#})  \\ && +\alpha(v)(u\times \alpha^{\#}) -
\alpha(u)(v\times \alpha^{\#}) +\varphi(u,v,\alpha^{\#})\alpha^{\#}\;]
\end{eqnarray*}

So the equation  $(u \times v)_{\lambda}(x)=(u\times v ) (x)$ is given by the equation $F=0$ where:
 \begin{eqnarray*}
F= a \chi(u,v,\alpha^{\#})  -\alpha(v)(u\times \alpha^{\#})+
\alpha(u)(v\times \alpha^{\#}) -\varphi(u,v,\alpha^{\#})\alpha^{\#}\ + |\alpha|^2(u\times v )
 \end{eqnarray*}

\vspace{.05in} Note that when $\alpha^{\#} \in {\bf E}$, the
equation $F(x)=0$ holds for all $x\in M$. Let us choose our
deformation $\alpha^{\#} \in {\bf V}$, which is a transverse
section of $ {\bf V}\to M$.  In this case by (4), Lemma 1, and
Lemma 2 the equation $F(x)=0 $ is equivalent  to $$a
J(\alpha^{\#})= -|\alpha|^{2}(u\times v)$$ where $J$ is the
complex structure defined in Lemma 2. Since $J(\alpha^{\#})\in
{\bf V}$ and $u\times v\in {\bf E}$, this equation holds only  at
points satisfying $\alpha^{\#}(x)=0$. By taking derivative of
$F(a,\alpha)$ we see that $F$ is transversal to
$\tilde{M}_{\varphi }$ when  $a\neq 0$. \qed

\vspace{.1in}

\section{ Cayley Submanifolds of $Spin(7)$}

 Much of what we have discussed  for  associative submanifolds of a $G_2$ manifold holds for Cayley submanifolds of a $Spin(7)$ manifold. Let  $(N^{8},\Psi)$ be a $Spin(7)$ manifold, and $\CP_{Spin(7)} (N)\to N$ be its $Spin(7)$  frame bundle, and $G(4,8)$ be the Grassmannian of oriented $4$ planes in $\R^8$. As in the $G_2$ case we can form the bundle
$$ \tilde{N}= \CP(N)\times_{SO(8)} G(4,8) \to N.$$
\noindent Similarly we have the universal bundles $\Xi $ , $\V \to
\tilde{N}$ which are  fiberwise extensions of the canonical bundle
$\xi \to G(4,8)$  and its dual $\nu=\xi^{\perp} \to G(4,8)$,
respectively. $ Hom (\Xi, \V)= \Xi^{*}\otimes \V \to \tilde{N}$ is
the vertical subbundle of $T(\tilde{N}) \to N$ with fibers
$TG(4,8)$. Let $G^{\Psi} (4,8)$ be the Grassmannian of Cayley
$4$-planes in $G(4,8)$ consisting of elements $L\in G(4,8)$
satisfying $\Psi|L=vol(L)$. The group  $Spin(7)$ acts transitively
on  $G^{\Psi} (4,8)$ with the stabilizer  $(SU(2) \times SU(2)
\times SU(2))/{\Z_2}$. Therefore, $G^{\Psi} (4,8)$ can be
identified by the quotient of $Spin(7)$ with the subgroup
$$ (SU(2)\times SU(2)\times SU(2))/\Z_{2} \subset Spin(7).$$
The action of $[q_{+},q_{-},\lambda ] \in (SU(2) \times SU(2) \times SU(2))/{\Z_2}$  on $\R^8=\H \oplus  \H$ is given by $(x,y)\to (q_{+}xq_{-}^{-1}, q_{+}y \lambda^{-1})$.
As in $G_2$ case there is  the  Cayley Grassmannian bundle
$$\tilde{N}_{\Psi} = \CP_{Spin(7)} (N)\times _{Spin(7)}G^{\varphi}(4,8) \to N$$
which is  $\tilde {N}_{\Psi }=\CP(N)/(SU(2) \times SU(2) \times SU(2))/{\Z_2}\to \CP(N)/Spin(7) =N$. We have restriction of the bundles
$\Xi^{*}$, $\V \to \tilde{N}_{\Psi}\subset \tilde{N}$.  Furthermore,  the principal $(SU(2) \times SU(2) \times SU(2))/{\Z_2}$  bundle $\CP(N)\to \tilde {N}_{\Psi }$ gives the following associated vector bundles over
$ \tilde {N}_{\Psi }$ via the  representations
(see  \cite{hl}, \cite{m}).
\begin{equation}
\begin{array}{lcc}
\W^{+} =\V :&y \mapsto  q_{+}y\lambda^{-1}   &  \hspace{1in} \\
\W^{-} :&y \mapsto  q_{-}y\lambda^{-1}   &  \hspace{1in} \\
\Xi^{*}:  & x\mapsto  q_{+}xq_{-}^{-1} &   \\
\lambda_{+}(\Xi^{*}) :& x\mapsto q_{+}xq_{+}^{-1}\\
   \lambda_{-}(\Xi^{*}): & x \mapsto q_{-}xq_{-}^{-1}&\\
  \lambda_{-}(\W^{+}):  \;  & x \mapsto  \lambda x  \lambda ^{-1} & \\
  \end{array}
 \end{equation}
where $[q_{+},q_{-},\lambda]\in SU(2)\times SU(2)\times SU(2))/ \Z_{2}$. We can identify:
$ \; \lambda_{+}(\W^{+})= \lambda_{+}(\Xi^{*})$,  and  we have the usual decomposition $\Lambda^{2}(\Xi^{*})=\lambda_{+}(\Xi^*)\oplus \lambda_{-}(\Xi^{*})$. 
We have the Clifford multiplications $\Xi^{*}\otimes \W^{\pm}\to\ \W^{\mp}$ given by:  $x\otimes y \mapsto -\bar{x}y$ and $x\otimes y \mapsto xy$,  on $\W^{+}$ and $\W^{-}$ respectively,  which extends to
$\Lambda^{2}(\Xi^{*})\otimes \W^{+}\to\ W^{+}$.

\vspace{.1in}

  The  Gauss map of an imbedding $f: X^{4}\hookrightarrow N^8$ of any $4$-manifold canonically lifts to an imbedding   $\tilde{f}: X^{4}\hookrightarrow \tilde{N}$, and the
pull backs $\tilde{f}^{*} \Xi^{*}=T^{*}(X)$ and  $\tilde{f}^{*}
\W^{+}=\nu(X)$ give cotangent and normal bundles of $X$.
Furthermore, if $X$ is a Cayley submanifold of N then the image of
$\tilde{f}$ lands in $ \tilde {N}_{\Psi }$;  in this case pulling
back the principal  $Spin(7)$ frame  bundle $\CP(N)\to N$ induces
an $(SU(2)\times SU(2)\times SU(2))/ \Z_{2}$ bundle $\CP (X)\to
X$. So by the representations (58) we get  associated vector
bundles $W^{+}=\nu(X), W^{-}, T^{*}(X)$ over $X$, i.e. the pull-backs of $\W^{+}, \W^{-}, \Xi^{*}$. So we have the actions $W^{+}\otimes
\lambda_{-}(W^{+})\to W^{+}$ and  $T^{*}X \otimes W^{\pm}\to
W^{\mp}$ and  $\Lambda^{2}(T^{*}X) \otimes W^{+}\to W^{+}$.

\vspace{.1in}

The  Levi-Civita connection induced by the $Spin(7)$ metric on
$N$,  induces  connections on tangent and normal bundle of any
submanifold $X^{4}\subset N$. Call these connections background
connections. Let $\A_{0}$ be the induced connection on
$\nu(X)=W^{+}$. Using the Lie algebra decomposition $
so(4)=so(3)\oplus so(3) $, we can decompose $ \A_{0}= S_{0}\oplus
A_{0}$, where $S_{0}$ and  $A_{0}$ are connections on
$\lambda_{+}(T^{*}X) $ and  $\lambda_{-}(W^{+})$, respectively.
Any connection $A$ of $ \lambda_{-}(W^{+})$  is in the form
$A=A_{0} +a$ where $a\in  \Omega^{1}(X,  \lambda_{-}(W^{+}))$, and
by  the association $ A\mapsto S_{0}\oplus A $  it  induces a
connection on $\nu(X)$. We will denote this connection   by
$\A=\A(a)$,  and   $\A_{0}=\A(0)$. Later we will consider
deformations
\begin{equation}
\A_{0} \mapsto \A.
\end{equation}

 \vspace{.1in}

Let  $\nabla_{\A}: \Omega^{0}(X, W^{+})\to \Omega^{1}(X, W^{+})$
by $ \nabla_{A}=\sum e^i \;\otimes \nabla_{e_i}$, where
$\{e_{i}\}$and $\{e^{i}\}$ are  orthonormal tangent and cotangent
frame fields of $X$, respectively. When $X$ is a Cayley manifold,
the Clifford multiplication  gives  the twisted Dirac operator:
 \begin{equation}
 \Di _{\A}: \Omega^{0}(X,W^{+}) \to \Omega^{0}(X,W^{-})
 \end{equation}
The  kernel of $\Di_{\A_{0}}$ gives the  infinitesimal deformations of Cayley submanifolds (\cite{m}). As in the associative case by deforming $\A_{0} \to \A$ we can make cokernel  of $\Di_{\A}$  zero.

  \vspace{.1in}

Similar to the  case of $\Lambda$-associative submanifolds in $G_2$
manifolds,  we can study the Cayley submanifolds in  $Spin(7)$
manifolds with complex normal bundles. There are several ways of
lifting various subbundles to complex bundles, for example
 $$  Spin^{c}(4)= ( SU(2)\times SU(2)\times S^1 )/ \Z_{2}\subset (SU(2)\times SU(2)\times SU(2))/\Z_{2}$$
 gives a $Spin^{c} (4)$ bundle $\CP(N)\to \bar{N}_{\Psi}=\CP(N)/Spin^{c}(4)$, and
 we have all the corresponding bundles  of (58) over  $\bar{N}_{\Psi}$ (except in this case we have $\lambda \in S^{1}$).  The $S^2$-bundle $\bar{N}_{\Psi }\to \tilde {N}_{\Psi }$ can  be identified with the sphere bundle  of  $ \lambda_{-}(\W^{+})\to \tilde{N}_{\Psi }$, and the sections of this bundle  correspond to almost complex structures on $\W^{\pm}$.   Previously, in the case of
$7$-manifolds,  existence of such sections  followed from the
existence of $2$-frame field  \cite{t}, in the $8$-dimensional
$Spin(7)$ case we don't have a clean analogue of \cite{t}, so in
this case we will make this an assumption and proceed. So consider
a $Spin(7)$ manifold $(N^8, \Psi, \Lambda)$ with a unit section
$\Lambda: \tilde{N}_{\Psi}\to \lambda_{-}(\W^{+})$.
 Hence $\W^{\pm}\to  \bar{N}_{\Psi}$ are $U(2)$ bundles, and $\lambda_{-}(W^+)$ is a line bundle $L\to  \bar{N}_{\Psi}$. As in (50) there is a  quadratic bundle map $\sigma:  \W^{+} \otimes {\W}^{+} \to \lambda_{+}(\Xi^{*})$   \begin{equation*}
 \sigma(x,x)=-\frac{1}{2}(x i \bar{x})i.
\end{equation*}
Now if  $f:X^4\hookrightarrow N^8$ is a  Cayley submanifold, we
can pull back these structures onto $X$ by  $\Lambda\circ
\tilde{f}$. Then we can ``perturb''  the local Cayley deformations
of $X$  by deforming the connection as in (59), i.e. the kernel of
the Dirac operator  of (60). Then if we  can  cut down the
solution space $\Di_{\A}^{-1}(0)$ by a second natural equation (by
using ``$a$'' as a free wariable) we arrive to  the Seiberg-Witten
equations:
\begin{equation}
  \begin{array}{ccc}
 \Di_{A}(v) &=& 0\\
 F^{+}_{A} &=  & \sigma (v,v)
\end{array}
 \end{equation}
As usual, by  perturbing these equations by  elements of
$\Omega^{2}_{+}(X)$,   i.e. by changing the second equation with
 $F^{+}_{A} +\delta = \sigma (v,v)$ with  $\delta \in \Omega^{2}_{+}(X)$  we get smoothness on the  zero locus of the parameterized equation $F=0$ where
 \begin{equation*}
F: \Omega^{0}(X,{W}^{+}) \times \CA(L) \times \Omega^{2}_{+} (X)\to \Omega^{0}(X,{W}^{-})\times \Omega^{2}_{+} (X)
\end{equation*}
 and by generic choice of $\delta$ we can make the solution set   $F_{\delta}^{-1}(0)$ smooth.  The normal bundle ${W}^{+}$ of $X$ has a complex structure, so  the  gauge group  $\CG(L)=Map (X, S^1)$ acts on   $F_{\delta}^{-1}(0) $, and makes quotient $F_{\delta}^{-1}(0) / \CG(L)$ a  smooth manifold whose dimension $d$ can be calculated from the index of the elliptic complex:
 \begin{equation}
\Omega^{0}(X) \to \Omega^{0}(X,{W}^{+}) \times \Omega^{1}(X)
\to \Omega^{0}(X,{W}^{-})\times \Omega^{2} _{+}(X)
\end{equation}
\noindent  where the first map comes from gauge group action. As  in Seiberg-Witten we get:
 \begin{equation}
d=\frac{1}{4}\;[ c_{1}^{2}(L)-(2 e(X) +3 \sigma(X)]
\end{equation}

\noindent Here $e$ and $\sigma$  denote the Euler characteristic
and the signature. In particular, these parametrized deformations of
complex Cayley submanifolds in $\Omega^{0}(X,W^{+})\times \CA(L)$
are unobstructed.

{\Thm Given  $(X, \Psi, \Lambda)$, to any Cayley submanifold $f:
X^4 \hookrightarrow N$ we can assign a number $\Lambda_{\Psi}
(X)\in  \Z$. Furthermore, the Seiberg-Witten equations of (61) can
be pulled back by $f$ from global equations on $N$ (analogue of
Proposition 9). }

\vspace{.1 in}

Note that $SU(3)$ and $G_{2}$ also act on the corresponding
special Lagrangian and coassociative Grassmannians with $SO(3)$
and $SO(4)$ stabilizers, respectively  \cite{hl}, giving the
identifications  $G^{SL}(3,6)=SU(3)/SO(3)$ and
$G^{coas}(4,7)=G_{2}/SO(4)$. As before, one can study special
Lagrangians in a Calabi-Yau manifold, and coassociative
submanifolds in a $G_2$ manifold, by lifting their normal bundles
to $SU(2)$. Their deformation spaces are unobstructed and can be
identified with $H^1$ and $H_{+}^{2}$, respectively. With a
similar approach we can relate them to the reduced Donaldson
invariants (as the $\Lambda$-associative and similarly defined Cayley's are related to Seiberg-Witten invariants).  Similarly one can treat
the deformations of associative submanifolds whose boundaries lie on coassociative submanifolds, and the Cayley's in $Spin(7)$ with
associative boundaries in $G_2$. Also asymptotically cylindrical
associative submanifolds in a $G_2$ manifold with a Calabi-Yau
boundary  have similar local deformation spaces, their
deformations are related to the corresponding holomorphic curves
inside the Calabi-Yau boundary.

\end{document}